\documentclass[12pt]{article}
\usepackage{amssymb}
\usepackage{amsmath}
\usepackage{amsfonts}

\newtheorem{theorem}{Theorem}

\newtheorem{example}[theorem]{Example}

\newtheorem{lemma}[theorem]{Lemma}

\newtheorem{proposition}[theorem]{Proposition}

\newenvironment{proof}[1][Proof]{\noindent\textbf{#1.} }{\ \rule{0.5em}{0.5em}}
\begin{document}

\title{Exponential and Gaussian behavior in the tails of multivariate
functions}
\author{Daniel J. Fresen\thanks{%
University of Pretoria, Department of Mathematics and Applied Mathematics,
daniel.fresen@up.ac.za or djfb6b@mail.missouri.edu MSC: 26B35, 60F05, 60G55,
60G70.}}
\date{}
\maketitle

\begin{abstract}
We observe that approximate copies of the function $\Lambda _{n}:\mathbb{R}%
^{n}\rightarrow (0,\infty )$ defined by%
\begin{equation*}
\Lambda _{n}(x)=\exp \left( -x_{1}-\pi \sum_{i=2}^{n}x_{i}^{2}\right) 
\end{equation*}%
appear in the tails of a large class of functions, with properties related
to coordinate independence, convexity, homotheticity, and homogeneity. The
function $\Lambda _{n}$ is an entropy maximizer (on a half-space) that is
uniquely determined by a homogeneity condition together with rotational
invariance about the $x_{1}$ direction and its behavior near the origin.

These results are connected to the limiting Poisson point processes found
near the edges of large random samples, as well as the conditioning of
random vectors on certain rare events, and can be thought of as variations
of Laplace's method for estimating integrals.
\end{abstract}

\section{Introduction}

An overview of the paper is as follows (we refer to Section \ref{main
results} for more details). Consider any $n\in \mathbb{N}$ with $n\geq 2$
and a function $f:\mathbb{R}^{n}\rightarrow (0,\infty )$. There are various
conditions that will be imposed on $f$ at different times that involve
smoothness, the growth of the derivatives of $-\log f$, convexity,
coordinate independence, homotheticity, and homogeneity. We are given a
point $y\in \mathbb{R}^{n}$ such that $|y|$ is large. Sometimes we also
insist that $y_{1}$ or $\min \{|y_{i}|:1\leq i\leq n\}$ is large. Consider
the function $\Lambda _{n}:\mathbb{R}^{n}\rightarrow (0,\infty )$ defined by%
\begin{equation}
\Lambda _{n}(x)=\exp \left( -x_{1}-\pi \sum_{i=2}^{n}x_{i}^{2}\right) 
\label{basic fnc def lambda}
\end{equation}%
We wish to find an invertible $n\times n$ matrix $T$ (depending on $y$ and $f
$) such that for all $x\in \mathbb{R}^{n}$,%
\begin{equation}
\frac{f\left( y+Tx\right) }{f(y)}=\left( 1+\delta _{y}(x)\right) \Lambda
_{n}(x)  \label{abcde random}
\end{equation}%
where $\delta _{y}(x)$ is an error term, usually given with appropriate
quantitative bounds. As it turns out, for a large class of functions it is
possible to achieve this approximation with $\delta _{y}(\cdot )\rightarrow 0
$ uniformly on compact subsets of $\mathbb{R}^{n}$, as $y\rightarrow \infty $
in the sense described above.

\begin{example}
Consider the function $f:\mathbb{R}^{2}\rightarrow (0,\infty )$ defined by $%
f(x)=e^{-(x_{1}^{4}+x_{2}^{4})}$, and the point $y=(t,t)$, for some $t\in 
\mathbb{R}$. Let%
\begin{equation*}
T=\left[ 
\begin{array}{cc}
a & -b \\ 
a & b%
\end{array}%
\right] 
\end{equation*}%
where $a=t^{-3}/8$ and $b=t^{-1}\sqrt{\pi /12}$. Then (after a little
algebra), $f(y+Tx)$ can be written as%
\begin{equation*}
\exp \left( -2t^{4}-x_{1}-\pi x_{2}^{2}-\frac{x_{1}^{4}}{2048t^{12}}-\frac{%
\pi x_{1}^{2}x_{2}^{2}}{64t^{8}}-\frac{\pi ^{2}x_{2}^{4}}{72t^{4}}-\frac{%
x_{1}^{3}}{64t^{8}}-\frac{\pi x_{1}x_{2}^{2}}{4t^{4}}-\frac{3x_{1}^{2}}{%
16t^{4}}\right) 
\end{equation*}%
As $t\rightarrow \infty $, the function $x\mapsto f(y+Tx)/f(y)$ converges
uniformly on compact subsets of $\mathbb{R}^{2}$ to the function $\exp
\left( -x_{1}-\pi x_{2}^{2}\right) $.
\end{example}

\subsection{\label{interpret MEVT}Interpretation: a multivariate extreme
value theory}

We refer the reader to \cite{BGST, Col, deHF, FaHuRe, Res} for background.
Classical extreme value theory in a single dimension is concerned with the
behavior of a large i.i.d. random sample $(X_{i})_{1}^{N}$ around its
maximum and minimum values. Let $F\left( x\right) =P\left\{ X_{1}\leq
x\right\} $. By the Fisher-Tippett-Gnedenko theorem, if a limiting
distribution for $\max (X_{i})_{1}^{N}$ exists as $N\rightarrow \infty $,
then it must be one of 3 fundamental types. More precisely, if there exist
sequences $(a_{N})_{1}^{\infty }\subset (0,\infty )$ and $%
(b_{N})_{1}^{\infty }\subset \mathbb{R}$ such that for all $x\in \mathbb{R}$,%
\begin{equation}
\lim_{N\rightarrow \infty }\mathbb{P}\left\{ \frac{\max (X_{i})_{1}^{N}-b_{N}%
}{a_{N}}\leq x\right\} =\lim_{N\rightarrow \infty }F\left(
a_{N}x+b_{N}\right) ^{N}=G(x)  \label{classical max}
\end{equation}%
and $G$ is non-degenerate (i.e. is not the C.D.F. of a Dirac point mass),
then the sequences $(a_{N})_{1}^{\infty },(b_{N})_{1}^{\infty }$ can be
chosen such that $G$ is of the form%
\begin{eqnarray*}
\Phi _{p}(x) &=&\left\{ 
\begin{array}{ccc}
0 & : & x<0 \\ 
\exp \left( -x^{-p}\right) & : & x\geq 0%
\end{array}%
\right. \hspace{0.35in}\Psi _{p}(x)=\left\{ 
\begin{array}{ccc}
\exp \left( -\left\vert x\right\vert ^{p}\right) & : & x<0 \\ 
1 & : & x\geq 0%
\end{array}%
\right. \\
\Lambda (x) &=&%
\begin{array}{ccc}
\exp \left( -e^{-x}\right) & : & x\in \mathbb{R}%
\end{array}%
\end{eqnarray*}%
where $p>0$. This can be expressed in the language of point processes. Let $%
\delta (x)$ denote the Dirac point mass at a point $x$ (in whatever space $x$
exists). The random measure%
\begin{equation}
\nu _{N}=\sum_{i=1}^{N}\delta \left( \frac{X_{i}-b_{N}}{a_{N}},\frac{i}{N}%
\right)  \label{random 2D measure}
\end{equation}%
converges (in an appropriate sense, see for example \cite[Definition 7.1,
Theorem 7.1]{Col} and \cite[Section 3.5, Proposition 3.21]{Res}) to a
non-homogeneous two dimensional Poisson point process $\nu ^{\ast }$ on $%
\left\{ x\in \mathbb{R}:0<G\left( x\right) <1\right\} \times \left(
0,1\right) $ with intensity%
\begin{equation*}
\mathbb{E}\nu ^{\ast }\left[ \left( u_{1},v_{1}\right] \times \left(
u_{2},v_{2}\right] \right] =\left( u_{2}-v_{2}\right) \int_{u_{1}}^{v_{1}}%
\frac{G^{\prime }\left( s\right) }{G\left( s\right) }ds
\end{equation*}%
Extreme value theory is intimately linked to the theory of regular variation
and $\Gamma $-variation. Going back to de-Haan, see e.g. \cite{deHa74, Res},
a function $f:\mathbb{R}\rightarrow \left( 0,\infty \right) $ is of class $%
\Gamma $ if there exists $t:\mathbb{R}\rightarrow \left( 0,\infty \right) $
such that for all $x\in \mathbb{R}$,%
\begin{equation}
\lim_{y\rightarrow \infty }\frac{f\left( y+t\left( y\right) x\right) }{%
f\left( y\right) }=e^{x}  \label{original gamma var}
\end{equation}%
Whether one uses $e^{x}$ or $e^{-x}$ is a fairly superficial matter.
Convergence of the re-scaled $\max $ to a limiting distribution as in (\ref%
{classical max}) can be expressed in terms of regular or $\Gamma $
variation. A limit exists in (\ref{classical max}) with $G=\Lambda $ if and
only if there exists $t:\mathbb{R}\rightarrow \left( 0,\infty \right) $ such
that for all $x\in \mathbb{R}$, 
\begin{equation*}
\lim_{y\rightarrow \infty }\frac{1-F\left( y+t\left( y\right) x\right) }{%
1-F\left( y\right) }=e^{-x}
\end{equation*}%
which essentially follows from%
\begin{equation*}
F\left( a_{N}x+b_{N}\right) ^{N}=\left[ 1-\left( 1-F\left( b_{N}\right)
\right) \frac{1-F\left( a_{N}x+b_{N}\right) }{1-F\left( b_{N}\right) }\right]
^{N}
\end{equation*}%
taking $b_{N}$ such that $N\left( 1-F\left( b_{N}\right) \right) \rightarrow
1$ and $a_{N}=t\left( b_{N}\right) $. This is implied by $X_{1}$ having a
density $f$ such that%
\begin{equation}
\lim_{y\rightarrow \infty }\frac{f\left( y+t\left( y\right) x\right) }{%
f\left( y\right) }=e^{-x}  \label{gamma variation of density}
\end{equation}%
To extend extreme value theory to a multivariate setting, the typical
approach that has been taken is to study the limiting distribution (after
re-scaling) of the coordinatewise maximum,%
\begin{equation*}
\max (X_{i})_{1}^{N}=\left( \max_{1\leq i\leq N}X_{i,j}\right) _{j=1}^{n}
\end{equation*}%
where $(X_{i})_{1}^{N}$ is a sequence of random vectors in $\mathbb{R}^{n}$,
and $X_{i}=(X_{i,j})_{j=1}^{n}$. Such a multivariate theory focuses on a
classification of possible limiting distributions, their marginals, and the
dependence structure that exists within these distributions. Other
interpretations of multivariate extreme value theory involve the convex
hull, see for example \cite[p. ix]{Res} and references in \cite{Fr0, FrVi},
as well as vertices of the convex hull and Pareto points, see \cite[Section
7.2]{FaHuRe}. The results of this paper can be interpreted as a multivariate
extension of extreme value theory that goes in a different direction, more
in line with (\ref{random 2D measure}), that we now discuss.

There are at least two advantages of studying the processes $\nu _{N}$ as
opposed to the re-scaled maximum $a_{N}^{-1}\left( \max
(X_{i})_{1}^{N}-b_{N}\right) $. The first is that the process $\nu _{N}$
(and the limiting process $\nu ^{\ast }$) give more explicit and readily
available information about other sample points $X_{i}<\max (X_{i})_{1}^{N}$%
. For example, for any fixed $k\in \mathbb{N}$,%
\begin{equation*}
\lim_{N\rightarrow \infty }\mathbb{P}\left\{ \frac{X_{(N-k+1)}-b_{N}}{a_{N}}%
\leq x\right\} =G\left( x\right) \sum_{j=0}^{k-1}\frac{\left( -\log G\left(
x\right) \right) ^{j}}{j!}
\end{equation*}%
where $X_{\left( N-k+1\right) }$ denotes the $k^{th}$ largest of the sample $%
\left( X_{i}\right) _{1}^{N}$. This comes down to convergence of the
probabilities%
\begin{equation*}
\lim_{N\rightarrow \infty }\mathbb{P}\left\{ \nu _{N}\left( x,\infty \right)
\leq k-1\right\} =\mathbb{P}\left\{ \nu ^{\ast }\left( x,\infty \right) \leq
k-1\right\}
\end{equation*}%
The joint distribution may be calculated similarly. While it is simple
enough to calculate the limiting distribution of $X_{(N-k+1)}$ from scratch,
the point is that this information is already contained in the behavior of $%
\nu _{N}$. The second advantage is that the processes $\nu _{N}$ make no
reference to an ordering. A multivariate theory based on limiting point
processes closely parallels the univariate theory: the only difference is
that the affine transformation $x\mapsto b_{N}+a_{N}x$ and its inverse $%
x\mapsto a_{N}^{-1}\left( x-b_{N}\right) $ are now replaced with the
transformation $x\mapsto y+Tx$ and its inverse $x\mapsto T^{-1}\left(
x-y\right) $.

\begin{example}[similar to existing results]
\label{reg var Poisson limit}Let $\mu $ be a probability measure and $\eta $
a non-negative Radon measure, both on the Borel subsets of $\mathbb{R}%
^{n}\backslash \left\{ 0\right\} $ (with its subspace topology inherited
from $\mathbb{R}^{n}$), such that $\lim_{t\rightarrow \infty }t^{p}\mu
\left( tE\right) =\eta \left( E\right) $ for all Borel sets $A\subset 
\mathbb{R}^{n}\backslash \left\{ 0\right\} $ with $\eta \left( \partial
A\right) =0$ (i.e. $t^{p}\mu \left( t\cdot \right) $ converges weakly to $%
\eta $) for some $p>0$, and let $\left( X_{i}\right) _{1}^{\infty }$ be an
i.i.d. sequence of random vectors, each with distribution $\mu $. The
sequence of processes%
\begin{equation*}
\nu _{N}=\sum_{i=1}^{N}\delta \left( N^{-1/p}X_{i}\right)
\end{equation*}%
then converges weakly (in the sense of \cite[Section 3.5]{Res}) to a Poisson
point process with intensity $\eta $. This follows from \cite[proof of
Proposition 3.21]{Res}, since $N\mathbb{P}\left\{ N^{-1/p}X_{1}\in \cdot
\right\} =N\mu \left( N^{1/p}\cdot \right) \rightarrow ^{w}\eta $ on $%
\mathbb{R}^{n}\backslash \left\{ 0\right\} $. This holds when, for example,%
\begin{equation*}
\frac{d\mu }{dx}=\prod_{i=1}^{n}\left( 1+\left\vert x_{i}\right\vert \right)
^{-1-p_{i}}\hspace{0.65in}\frac{d\eta }{dx}=\prod_{i=1}^{n}\left\vert
x_{i}\right\vert ^{-1-p_{i}}
\end{equation*}%
where $p_{i}>0$ and $\sum p_{i}=p$. See related results in \cite[Section 2]%
{ColTaw} \cite[Eq. (5.4.15)]{Miko04} \cite[Theorem 2]{Foug04} and references
therein.
\end{example}

An interpretation of our work is that for large $N$ the random measure%
\begin{equation*}
\nu _{N}=\sum_{i=1}^{N}\delta \left( T^{-1}(X_{i}-y)\right)
\end{equation*}%
where $\left( X_{i}\right) _{1}^{N}$ is now an i.i.d. sequence of random
vectors, approximates a Poisson point process with (density of) intensity $%
\exp \left( -x_{1}-\pi \sum_{i=2}^{n}x_{i}^{2}\right) $. Here $y\in \mathbb{R%
}^{n}$ is a point near the edge of the random sample (the same $y$ as
before), and $T\in GL(n)$ depends on $y$, and we assume $f(y)\approx \left(
N\det (T)\right) ^{-1}$. The last condition is what determines how far out
in the tails $y$ should be. We assume that the common distribution of each $%
X_{i}$ has a density $f=d\mu /dx$, which is the function to which we apply
the results of the paper. This interpretation can be expressed more
precisely as follows.

\begin{theorem}[Poissonization]
\label{extremal proc} Let $n\geq 2$ and let $\mu $ be a probability measure
on $\mathbb{R}^{n}$ with density $f=d\mu /dx$. There exists a probability
space $(\Omega ,\Im ,\mathbb{P})$ and an i.i.d. sequence of random vectors $%
(X_{i})_{1}^{\infty }$ in $\mathbb{R}^{n}$ defined on $\Omega $ with common
distribution $\mu $ such that the following holds. Let $N\in \mathbb{N}$
with $N>N_{0}$ ($N_{0}$ a universal constant), $y\in \mathbb{R}^{n}$, $T\in
GL(n)$, $0<\varepsilon <1$, and $R>0$, and consider the random measure%
\begin{equation*}
\nu _{N}=\sum_{i=1}^{N}\delta \left( T^{-1}(X_{i}-y)\right)
\end{equation*}%
Let $B^{\circ }(0,R)=\left\{ x\in \mathbb{R}^{n}:\left\vert x\right\vert
<R\right\} $. Assume that%
\begin{equation}
\int_{\left\vert x\right\vert <R}\left\vert \exp \left( -x_{1}-\pi
\sum_{2}^{n}x_{i}^{2}\right) -N\det (T)f\left( y+Tx\right) \right\vert
dx<\varepsilon  \label{tot var approx}
\end{equation}%
Then there exists a Poisson point process $\nu =\nu _{\omega }$ ($\omega \in
\Omega $) on $B^{\circ }(0,R)$ with (density of) intensity $\exp \left(
-x_{1}-\pi \sum_{2}^{n}x_{i}^{2}\right) $ such that with probability at least%
\begin{equation*}
1-\varepsilon -2N^{-1/4}-3\sqrt{\frac{\log N}{N}}\left( 1+\int_{\left\vert
x\right\vert <R}\exp \left( -x_{1}-\pi \sum_{2}^{n}x_{i}^{2}\right) dx\right)
\end{equation*}%
for all Borel sets $E\subseteq B^{\circ }(0,R)$,%
\begin{equation*}
\nu _{N}(E)=\nu (E)
\end{equation*}
\end{theorem}

Note that the type of approximation one gets, $\nu _{N}|_{B^{\circ }\left(
0,R\right) }=\nu |_{B^{\circ }\left( 0,R\right) }$ with high probability, is
beyond what is possible in the setting of weak convergence (where $\nu
_{N}|_{B^{\circ }\left( 0,R\right) }$ and $\nu |_{B^{\circ }\left(
0,R\right) }$ often have disjoint support). This is facilitated through the
assumption (\ref{tot var approx}) on the total variation distance. Another
notable difference between our results (interpreted through Theorem \ref%
{extremal proc}) and results in the spirit of Example \ref{reg var Poisson
limit} is as follows: in Example \ref{reg var Poisson limit} the observer is
stationed at $y=0$ and we get a single global view of the entire sample,
whereas in our results the observer is stationed in the tails of the
probability distribution, near the outer regions of the point cloud $\left\{
X_{i}\right\} _{1}^{N}$, and we view the cloud locally.

Since the relationship between a large i.i.d. sample of fixed size and a
Poisson process is well understood, we don't include any further discussion
about limiting point processes and instead focus entirely on the analysis of
the density $f$. Obviously, (\ref{abcde random}) can be thought of as a
multivariate equivalent of de Haan's $\Gamma $-variation, as in (\ref%
{original gamma var}) and (\ref{gamma variation of density}). In this
context, both%
\begin{equation*}
\Lambda _{n}^{\ast }\left( x\right) =\exp \left( -x_{1}\right) \hspace{0.65in%
}\Lambda _{n}^{\sharp }\left( x\right) =\exp \left(
-\sum_{i=1}^{n}x_{i}\right)
\end{equation*}%
could also be thought of as multivariate analogues of $x\mapsto \exp \left(
-x\right) $, $x\in \mathbb{R}$. These two functions are equivalent up to
composition on the right with a linear map and are degenerate versions of $%
\Lambda _{n}$. The latter has a richer shape since it contains approximate
copies of $\Lambda _{n}^{\ast }$ but not vice versa (although this means
that copies of the function $\Lambda _{n}^{\ast }$ are slightly more
ubiquitous), and it models the behavior of a strictly convex point cloud.
Furthermore, the scale on which we view the point cloud in Theorem \ref%
{extremal proc} is approximately uniquely determined (up to orthogonal
transformations of $x_{2}\ldots x_{n}$), but this is not the case if we use $%
\Lambda _{n}^{\ast }$ with $n\geq 3$ since it is invariant under arbitrary
transformations of $x_{2}\ldots x_{n}$, in particular by elements of $%
SL\left( n-1\right) $.\ See \cite[Section 5.4.2]{Res} and \cite[Section 8.4]%
{BGST} for information on multivariate regular variation and references
contained therein (regular variation is the polynomial decay version of $%
\Gamma $-variation).

\subsection{\label{interp rare}Interpretation: conditioning on rare events}

Let $X$ be a single random vector in $\mathbb{R}^{n}$ with coordinates $%
\left( X_{i}\right) _{1}^{n}$ (a change in notation from Section \ref%
{interpret MEVT}) and joint density function $f$, and let $\theta \in
S^{n-1} $. Under appropriate assumptions on $f$ and $\theta $ (that depend
on which of our results is being invoked), if $L>0$ is large and we let $%
\widetilde{X} $ be the random vector $X$ conditioned on the event $\left\{
\sum_{1}^{n}\theta _{i}X_{i}\geq L\right\} $, then we may write%
\begin{equation*}
\widetilde{X}=y+TW
\end{equation*}%
where $y$ maximizes $f$ on $\left\{ x\in \mathbb{R}^{n}:\sum_{1}^{n}\theta
_{i}x_{i}=L\right\} $ and $W$ is a random vector in $\left[ 0,\infty \right)
\times \mathbb{R}^{n-1}$ with a density $h$ such that $\left\vert h\left(
x\right) -\exp \left( -x_{1}-\pi \sum_{i=2}^{n}x_{i}^{2}\right) \right\vert
\leq \varepsilon _{f,L,\theta }\left( x\right) $ for all $x\in \left[
0,\infty \right) \times \mathbb{R}^{n-1}$ and $\varepsilon _{f,L,\theta
}\left( \cdot \right) \rightarrow 0$ uniformly on compact sets as $%
L\rightarrow \infty $ and $\theta $ remains fixed (i.e. $W$ is approximately
a random vector with independent exponential and normal coordinates). If we
apply this interpretation to Theorem \ref{th log concave product}, for
example, then we must impose the condition that $\theta _{i}\neq 0$ for all $%
i$, which is not much of an assumption since if some $\theta _{i}=0$ then we
just drop that coordinate entirely from the problem (in this sense we
ultimately need two nonzero coordinates). As $L$ and $y$ are fixed and $%
x\rightarrow \infty $, the bound $\varepsilon _{f,L,\theta }\left( x\right) $
may become large, and we must impose some regularity on $f$ to prevent bad
behavior of $f\left( y+Tx\right) $. For example, if we assume that $f$ has a
log-concave density, i.e. $f=\exp \left( -g\right) $ for some convex
function $g$, then $h$ too will be log-concave. If $\varepsilon _{f,L,\theta
}\left( x\right) \rightarrow 0$ uniformly in $\left\{ x\in \left[ 0,\infty
\right) \times \mathbb{R}^{n-1}:\left\vert x\right\vert <R\right\} $, for
some $R>0$ with $R\rightarrow \infty $ as $L\rightarrow \infty $, then $h$
must continue to decay rapidly outside $\left\{ x\in \left[ 0,\infty \right)
\times \mathbb{R}^{n-1}:\left\vert x\right\vert <R\right\} $, and so%
\begin{equation*}
\lim_{L\rightarrow \infty }\frac{\int_{\left\{ x\in \left[ 0,\infty \right)
\times \mathbb{R}^{n-1}:\left\vert x\right\vert \geq R\right\} }h\left(
x\right) dx}{\int_{\left\{ x\in \left[ 0,\infty \right) \times \mathbb{R}%
^{n-1}:\left\vert x\right\vert <R\right\} }h\left( x\right) dx}=0
\end{equation*}%
This is important, otherwise the mass outside $\left\{ x\in \left[ 0,\infty
\right) \times \mathbb{R}^{n-1}:\left\vert x\right\vert <R\right\} $ may
dominate the conditioning that defines $\widetilde{X}$.

\begin{example}
\label{eg first condition}Let $X$ be a random vector in $\mathbb{R}^{2}$
with independent coordinates, each with distribution $\mathbb{P}\left\{
X_{i}\leq t\right\} =\left( 2\Gamma \left( 5/4\right) \right)
^{-1}\int_{-\infty }^{t}\exp \left( -u^{4}\right) du$. The joint density of $%
X$ is therefore%
\begin{equation*}
f\left( x_{1},x_{2}\right) =\left( 2\Gamma \left( 5/4\right) \right)
^{-2}\exp \left( -x_{1}^{4}-x_{2}^{4}\right)
\end{equation*}%
Let us condition on the event that $X_{1}+2X_{2}\geq L$, and label the
resulting random vector $\widetilde{X}$. The density $f$ is maximized on $%
\left\{ \left( x_{1},x_{2}\right) :x_{1}+2x_{2}=L\right\} $ at the point $%
y=\left( y_{1},y_{2}\right) $ with%
\begin{equation*}
y_{1}=\left( 1+2^{4/3}\right) ^{-1}L\hspace{0.65in}y_{2}=2^{1/3}\left(
1+2^{4/3}\right) ^{-1}L
\end{equation*}%
and setting $g=-\log f$, $\nabla g\left( y_{1},y_{2}\right) =4\left(
1+2^{4/3}\right) ^{-3}L^{3}\left( 1,2\right) $. Referring to the statement
and proof of Theorem \ref{th log concave product}, the norm $\left\Vert
\cdot \right\Vert _{\sharp }$ becomes%
\begin{equation*}
\left\Vert u\right\Vert _{\sharp }=\sqrt{\frac{6}{\pi }\left(
y_{1}^{2}u_{1}^{2}+y_{2}^{2}u_{2}^{2}\right) }
\end{equation*}%
and we may take%
\begin{equation*}
Q\left( x_{2}\right) =x_{2}\sqrt{\frac{\pi }{6\left(
4y_{1}^{2}+y_{2}^{2}\right) }}\left( 2,-1\right)
\end{equation*}%
The linear map $T$ then takes the form%
\begin{equation*}
Tx=\frac{1}{20}\left( 1+2^{4/3}\right) ^{3}L^{-3}x_{1}\left( 1,2\right)
+\left( 1+2^{4/3}\right) L^{-1}x_{2}\sqrt{\frac{\pi }{6\left(
4+2^{2/3}\right) }}\left( 2,-1\right)
\end{equation*}%
and $\widetilde{X}=\left( \widetilde{X}_{1},\widetilde{X}_{2}\right) $ where%
\begin{eqnarray*}
\widetilde{X}_{1} &=&\left( 1+2^{4/3}\right) ^{-1}L+\frac{1}{20}\left(
1+2^{4/3}\right) ^{3}L^{-3}Y+\frac{\left( 1+2^{4/3}\right) L^{-1}}{\sqrt{%
3\left( 4+2^{2/3}\right) }}Z \\
\widetilde{X}_{2} &=&2^{1/3}\left( 1+2^{4/3}\right) ^{-1}L+\frac{1}{10}%
\left( 1+2^{4/3}\right) ^{3}L^{-3}Y-\frac{\left( 1+2^{4/3}\right) L^{-1}}{2%
\sqrt{3\left( 4+2^{2/3}\right) }}Z
\end{eqnarray*}%
and $\left( Y,Z\right) $ converges weakly and in the total variation
distance (as $L\rightarrow \infty $) to a random vector in $\left[ 0,\infty
\right) \times \mathbb{R}$ with independent coordinates, the first standard
exponential and the second standard normal.
\end{example}

Example \ref{eg first condition} demonstrates a significant difference
between the univariate setting and the multivariate setting. In the former,
conditioning on $\left\{ X_{1}\geq T\right\} $, there is no normal
component, only exponential. In the latter, there are exponential and normal
components, and in this particular example the normal components dominate.
Conditioning on certain rare events thus facilitates convergence to a normal
even with $n=2$, unlike the classical central limit theorems where $n$ must
be large.

\subsection{\label{interpret Lapl}Interpretation: variations of Laplace's
method}

One reason for considering the limiting shape of functions in their tails is
that it leads to accurate estimates for $\int_{E}f(z)dz$ when $E\subset 
\mathbb{R}^{n}$ is such that (\ref{abcde random}) holds for all $x\in
T^{-1}\left( E-y\right) $. Assuming $\mathrm{vol}_{n}(E)>0$, such estimates
take the form%
\begin{equation}
\int_{E}f(z)dz=(1+\varepsilon )f(y)\left\vert \det (T)\right\vert
\int_{T^{-1}(E-y)}\exp \left( -x_{1}-\pi \sum_{i=2}^{n}x_{i}^{2}\right) dx
\label{tail integral}
\end{equation}%
where $y$ depends on $f$ and $E$, and $\varepsilon $ depends on $\delta
_{y}\left( \cdot \right) $ from (\ref{abcde random}). Similar estimates also
hold when the dimension of $E$ is less than $n$, for example when $E$ is a
hyperplane. An obvious application of (\ref{tail integral}) is to finding
asymptotic $\left( 1+o\left( 1\right) \right) $ estimates for tail
probabilities and densities of linear and nonlinear functionals of random
vectors, as in Example \ref{nonlinear deviation} below.

This is precisely the setting of Laplace's method, where $f:\mathbb{R}%
^{n}\rightarrow \mathbb{R}$ is of the form%
\begin{equation*}
f(t,z)=h(z)\exp \left( -tq(z)\right)
\end{equation*}%
for $t\in \mathbb{R}$, $z\in \mathbb{R}^{n-1}$, and we identify $\mathbb{R}%
\times \mathbb{R}^{n-1}$ with $\mathbb{R}^{n}$. Let us assume for simplicity
that $h$ is bounded and continuous, that $q$ is $C^{2}$\ and achieves a
strict global minimum at a point $z_{0}\in \mathbb{R}^{n-1}$, that the
Hessian $H_{q}(z_{0})$ is positive definite, that $q(z)\rightarrow \infty $
as $\left\vert z\right\vert \rightarrow \infty $, $\int_{\mathbb{R}^{n}}\exp
\left( -q(z)\right) dz<\infty $, and that $h(z_{0})\neq 0$. The objects $y$
and $T$ then take the form%
\begin{equation*}
y=(t,z_{0})\hspace{0.65in}T(s,u)=\left( \frac{s}{q(z_{0})},\sqrt{\frac{2\pi 
}{t}}H_{q}^{-1/2}(z_{0})u\right)
\end{equation*}%
where $s\in \mathbb{R}$ and $u\in \mathbb{R}^{n-1}$. The estimate (\ref%
{abcde random}) then leads to the classical approximation (setting $m=n-1$) 
\begin{equation*}
\int_{z\in \mathbb{R}^{m}}h(z)\exp \left( -tq(z)\right) dz=\int_{\mathcal{H}%
_{t}}f(x)dx=(1+o(1))\frac{h(z_{0})\exp \left( -tq(z_{0})\right) }{\sqrt{\det
H_{q}(z_{0})}}\left( \frac{2\pi }{t}\right) ^{m/2}
\end{equation*}%
as $t\rightarrow \infty $, where $\mathcal{H}_{t}=\left\{ x\in \mathbb{R}%
^{n}:x_{1}=t\right\} $ and integration is performed with respect to $(n-1)$
dimensional Lebesgue measure on $\mathcal{H}_{t}$.

\begin{example}[tail probabilities of a nonlinear functional]
\label{nonlinear deviation}Let $X$ be a random vector in $\mathbb{R}^{n}$
with distribution $\mu $ and density%
\begin{equation*}
f\left( x\right) =d\mu /dx=\left( 2\Gamma \left( 5/3\right) \right)
^{-n}\exp \left( -\sum_{i=1}^{n}\left\vert x_{i}\right\vert ^{3/2}\right)
\end{equation*}%
We estimate the tail probabilities $\mathbb{P}\left\{
\prod_{i=1}^{n}X_{i}\geq t\right\} $ in (\ref{nonlin tail prob}) below as $n$
is fixed and $t\rightarrow \infty $ (obviously in this particular example
one can convert to a linear functional by using the logarithm; the technique
is nonetheless demonstrated). For all $s>0$ let $\mathcal{M}_{s}=\left\{
x\in \left( 0,\infty \right) ^{n}:\prod_{i=1}^{n}x_{i}=s^{n}\right\} $. Note
that $r\mathcal{M}_{s}=\mathcal{M}_{rs}$ for all $r>0$. By log-concavity, $f$
is maximized over $\mathcal{M}_{t^{1/n}}$ at the point $y=\left(
t^{1/n},t^{1/n}\ldots ,t^{1/n}\right) $. Let $T\in GL\left( n\right) $ be as
in Theorem \ref{th log concave product}, in which case%
\begin{equation*}
\left\vert \mathrm{Det}T\right\vert =\frac{\left( 2\pi \right) ^{(n-1)/2}}{%
\sqrt{3n}}\left( \frac{4}{3}\right) ^{n/2}t^{-3/(4n)+1/4}
\end{equation*}%
As $t\rightarrow \infty $, the curvature of $\mathcal{M}_{t^{1/n}}$ at $y$
decreases to zero (at an appropriate rate), and the Hessian matrix of $%
g=-\log f$ at $y$ is the diagonal matrix with entries $\left( 3/4\right)
t^{-1/(2n)}$, so in operator norm%
\begin{equation*}
\left\Vert H_{g}\left( y\right) ^{-1/2}\right\Vert _{\ell
_{2}^{n}\rightarrow \ell _{2}^{n}}=\frac{2}{\sqrt{3}}t^{1/(4n)}\hspace{0.65in%
}\left\Vert T\right\Vert _{\ell _{2}^{n}\rightarrow \ell _{2}^{n}}=o\left(
t^{1/n}\right)
\end{equation*}%
This means that on the scale that matters, the curvature of $\mathcal{M}%
_{t^{1/n}}$ around $y$ becomes negligible. More precisely, the set $%
T^{-1}\left( E-y\right) $ increases to $\left( 0,\infty \right) \times 
\mathbb{R}^{n-1}$, where $E=\left\{ x\in \left( 0,\infty \right)
^{n}:\prod_{i=1}^{n}x_{i}\geq t\right\} $. By (\ref{tail integral}), $%
\mathbb{P}\left\{ \prod_{i=1}^{n}X_{i}\geq t\right\} =2^{n-1}\mu \left(
E\right) $ can be expressed as%
\begin{equation}
\frac{\left( 1+o_{n}\left( 1\right) \right) }{\sqrt{3n}}\left( 8\pi \right)
^{(n-1)/2}\left( \sqrt{3}\Gamma \left( 5/3\right) \right)
^{-n}t^{-3/(4n)+1/4}\exp \left( -nt^{3/(2n)}\right)  \label{nonlin tail prob}
\end{equation}%
as $t\rightarrow \infty $, where the subscript in $o_{n}\left( 1\right) $
indicates dependence on $n$.
\end{example}

\subsection{Further discussion}

Returning to the general setting, set $g=-\log f$ and assume momentarily
that $g$ is convex (in which case $f$ is said to be log-concave) and $C^{2}$%
. Let $H_{g}(y)$ denote the Hessian of $g$ at $y$. Provided $\nabla
g(y)^{\bot }\cap \mathrm{null}H_{g}(y)=\left\{ 0\right\} $, $f$ resembles
the function $\Lambda _{n}$ in a region surrounding $y$ when viewed in the
coordinate structure corresponding to the inner product%
\begin{equation*}
\left\langle u,v\right\rangle _{\sharp }=\frac{1}{2\pi }\left\langle
u,H_{g}(y)v\right\rangle +\left\langle u,\nabla g(y)\right\rangle
\left\langle v,\nabla g(y)\right\rangle
\end{equation*}%
More precisely, (\ref{abcde random}) holds with appropriate quantitative
bounds when $T\in GL\left( n\right) $ is such that for all $u,v\in \mathbb{R}%
^{n}$, $\left\langle Tu,Tv\right\rangle _{\sharp }=\left\langle
u,v\right\rangle $, and such that $Te_{1}=A^{-1}\nabla g(y)$, where%
\begin{equation*}
A=\frac{1}{2\pi }H_{g}(y)+\nabla g(y)\otimes \nabla g(y)\hspace{0.65in}%
A_{i,j}=\frac{1}{2\pi }\partial _{i,j}g(y)+\partial _{i}g(y)\partial _{j}g(y)
\end{equation*}%
In Theorem \ref{find scale} we provide a formula for the matrix $T$,
provided $f$ is smooth enough (in practice however, we don't always use this
exact $T$). Observe that $\nabla g(y)^{\bot }$ and $\mathrm{null}H_{g}(y)$
can be expressed using the following coordinate-free representations:%
\begin{eqnarray}
\nabla g(y)^{\bot } &=&\left\{ x\in \mathbb{R}^{n}:\left. \frac{d}{dt}%
g\left( y+tx\right) \right\vert _{t=0}=0\right\}  \notag \\
\mathrm{null}H_{g}(y) &=&\left\{ x\in \mathbb{R}^{n}:\left. \frac{d^{2}}{%
dt^{2}}g\left( y+tx\right) \right\vert _{t=0}=0\right\}
\label{loch hess monster}
\end{eqnarray}%
and that these subspaces are orthogonal with respect to $\left\langle \cdot
,\cdot \right\rangle _{\sharp }$. The case when $g$ is not convex is a
little different. For example setting $g(x_{1},x_{2})=x_{1}^{2}-x_{2}{}^{2}$%
, (\ref{loch hess monster}) clearly breaks down.

The function $\Lambda _{n}$ has, in a sense, the most natural possible shape
for a function that grows or decays rapidly. It has the homogeneity property
that for all $y\in \mathbb{R}^{n}$ there exists $T\in GL\left( n\right) $
such that for all $x\in \mathbb{R}^{n}$,%
\begin{equation*}
\frac{\Lambda _{n}\left( y+Tx\right) }{\Lambda _{n}(y)}=\Lambda _{n}(x)
\end{equation*}%
Furthermore $\Lambda _{n}$ is uniquely determined by this property, as well
as its gradient and Hessian at zero, $C^{3}$-smoothness, and the fact that $%
\Lambda _{n}\circ U=\Lambda _{n}$ for all $U\in O(n)$ such that $Ue_{1}=e_{1}
$ (to shorten the paper we have omitted the proof). It is also worth noting
that the function $\Lambda _{n}$ restricted to $[0,\infty )\times \mathbb{R}%
^{n-1}$ (which still encodes the full shape of $\Lambda _{n}$) maximizes
differential entropy among probability density functions of random vectors $X
$ supported on $[0,\infty )\times \mathbb{R}^{n-1}$ with mean $(1,0,\ldots 0)
$ and covariance%
\begin{equation*}
\left[ 
\begin{array}{cccc}
1 & 0 & \ldots  & 0 \\ 
0 & \left( 2\pi \right) ^{-1} & \ldots  & 0 \\ 
\vdots  & \vdots  & \ddots  & \vdots  \\ 
0 & 0 & \ldots  & \left( 2\pi \right) ^{-1}%
\end{array}%
\right] 
\end{equation*}

\section{\label{main results}Main results}

Recall that the density of an absolutely continuous $n$-fold product measure
on $\mathbb{R}^{n}$ can always be written in the form $f(x)=\prod_{i=1}^{n}%
\exp \left( -g_{i}(x_{i})\right) $, where each $g_{i}:\mathbb{R}\rightarrow
(-\infty ,\infty ]$.

\begin{theorem}[coordinate independence]
\label{th log concave product}Let $n\geq 2$ and for each $1\leq i\leq n$ let 
$g_{i}:\mathbb{R}\rightarrow \mathbb{R}$ be three times differentiable. Let $%
y\in \mathbb{R}^{n}$ such that for all $i$, $g_{i}^{\prime }(y_{i})\neq 0$
and $g_{i}^{\prime \prime }(y_{i})>0$, and consider the function $f:\mathbb{R%
}^{n}\rightarrow (0,\infty )$ defined by%
\begin{equation*}
f(x)=\prod_{i=1}^{n}\exp \left( -g_{i}(x_{i})\right)
\end{equation*}%
Set $g=-\log f$. There exists a linear bijection $T:\mathbb{R}%
^{n}\rightarrow \mathbb{R}^{n}$ such that for all $x\in \mathbb{R}^{n}$,%
\begin{equation}
\frac{f(y+Tx)}{f(y)}=\exp \left( -x_{1}-\pi \sum_{i=2}^{n}x_{i}^{2}+\delta
(x)\right)  \label{logconc prod bound}
\end{equation}%
where $\left\vert \delta (x)\right\vert $ can be bounded above by%
\begin{eqnarray*}
&&\frac{1}{2}\left\vert x_{1}\right\vert ^{2}\sum_{i=1}^{n}\frac{%
g_{i}^{\prime }(y_{i})^{4}}{\left\vert \nabla g(y)\right\vert ^{4}}\frac{%
g_{i}^{\prime \prime }(y_{i})}{g_{i}^{\prime }(y_{i})^{2}}+\sqrt{2\pi }%
\left\vert x_{1}\right\vert \left( \sum_{i=2}^{n}x_{i}^{2}\right)
^{1/2}\left( \sum_{i=1}^{n}\frac{g_{i}^{\prime }(y_{i})^{4}}{\left\vert
\nabla g(y)\right\vert ^{4}}\frac{g_{i}^{\prime \prime }(y_{i})}{%
g_{i}^{\prime }(y_{i})^{2}}\right) ^{1/2} \\
&&+\frac{1}{6}\left\vert x_{1}\right\vert ^{3}\sum_{i=1}^{n}\frac{%
g_{i}^{\prime }(y_{i})^{6}}{\left\vert \nabla g(y)\right\vert ^{6}}\frac{%
g_{i}^{\prime \prime }(y_{i})^{3/2}}{\left\vert g_{i}^{\prime
}(y_{i})\right\vert ^{3}}\frac{\left\vert g_{i}^{\prime \prime \prime
}(y_{i}+w_{i})\right\vert }{g_{i}^{\prime \prime }(y_{i})^{3/2}} \\
&&+\sqrt{\frac{\pi }{2}}\left\vert x_{1}\right\vert ^{2}\left(
\sum_{i=2}^{n}x_{i}^{2}\right) ^{1/2}\left( \sum_{i=1}^{n}\frac{%
g_{i}^{\prime }(y_{i})^{8}}{\left\vert \nabla g(y)\right\vert ^{8}}\frac{%
g_{i}^{\prime \prime }(y_{i})^{2}}{g_{i}^{\prime }(y_{i})^{4}}\frac{%
\left\vert g_{i}^{\prime \prime \prime }(y_{i}+w_{i})\right\vert ^{2}}{%
g_{i}^{\prime \prime }(y_{i})^{3}}\right) ^{1/2} \\
&&+\pi \left\vert x_{1}\right\vert \left( \sum_{i=2}^{n}x_{i}^{2}\right)
\max_{1\leq i\leq n}\left\{ \frac{g_{i}^{\prime }(y_{i})^{2}}{\left\vert
\nabla g(y)\right\vert ^{2}}\frac{g_{i}^{\prime \prime }(y_{i})^{1/2}}{%
\left\vert g_{i}^{\prime }(y_{i})\right\vert }\frac{\left\vert g_{i}^{\prime
\prime \prime }(y_{i}+w_{i})\right\vert }{g_{i}^{\prime \prime }(y_{i})^{3/2}%
}\right\} \\
&&+\frac{\sqrt{2\pi ^{3}}}{3}\left( \sum_{i=2}^{n}x_{i}^{2}\right)
^{3/2}\max_{1\leq i\leq n}\left\{ \frac{\left\vert g_{i}^{\prime \prime
\prime }(y_{i}+w_{i})\right\vert }{g_{i}^{\prime \prime }(y_{i})^{3/2}}%
\right\}
\end{eqnarray*}%
and $w\in \mathbb{R}^{n}$ is such that%
\begin{equation*}
\left( \sum_{i=1}^{n}g_{i}^{\prime \prime }(y_{i})w_{i}^{2}\right)
^{1/2}\leq \left\vert x_{1}\right\vert \left( \sum_{i=1}^{n}\frac{%
g_{i}^{\prime }(y_{i})^{4}}{\left\vert \nabla g(y)\right\vert ^{4}}\frac{%
g_{i}^{\prime \prime }(y_{i})}{g_{i}^{\prime }(y_{i})^{2}}\right) ^{1/2}+%
\sqrt{2\pi }\left( \sum_{i=2}^{n}x_{i}^{2}\right) ^{1/2}
\end{equation*}%
Furthermore, $T$ is of the form%
\begin{equation*}
Tz=\frac{\nabla g(y)}{|\nabla g(y)|^{2}}z_{1}+Q(z_{2},z_{3},\ldots ,z_{n})
\end{equation*}%
where $Q:\mathbb{R}^{n-1}\rightarrow \nabla f(y)^{\bot }$ and%
\begin{equation*}
\left\vert \mathrm{Det}T\right\vert =\left( 2\pi \right) ^{(n-1)/2}\left(
\sum_{i=1}^{n}\frac{g_{i}^{\prime }(y_{i})^{2}}{g_{i}^{\prime \prime }(y_{i})%
}\right) ^{-1/2}\left( \prod_{i=1}^{n}g_{i}^{\prime \prime }(y_{i})\right)
^{-1/2}
\end{equation*}
\end{theorem}

It is somewhat typical that in Theorem \ref{th log concave product} the
error term $\delta (\cdot )\rightarrow 0$ uniformly on compact subsets of $%
\mathbb{R}^{n}$ as $\min \left\{ \left\vert y_{i}\right\vert \right\}
\rightarrow \infty $. As Proposition \ref{pploy} below suggests, we can
often expect terms such as%
\begin{equation*}
\frac{g_{i}^{\prime \prime }(y_{i})}{g_{i}^{\prime }(y_{i})^{2}}\hspace{%
0.65in}\frac{\left\vert g_{i}^{\prime \prime \prime
}(y_{i}+w_{i})\right\vert }{g_{i}^{\prime \prime }(y_{i})^{3/2}}
\end{equation*}%
to be small when $\left\vert y_{i}\right\vert $ is large. This may then be
combined with the fact that%
\begin{equation*}
\sum_{i=1}^{n}\frac{g_{i}^{\prime }(y_{i})^{2}}{\left\vert \nabla
g(y)\right\vert ^{2}}=1
\end{equation*}

\begin{proposition}
\label{pploy}Let $\varepsilon >0$ and let $\omega :(0,\infty )\rightarrow 
\mathbb{R}$ be any differentiable function such that $\lim_{t\rightarrow
\infty }$ $t^{1/\varepsilon }\omega (t)=\infty $. Then%
\begin{equation}
\underset{t\rightarrow \infty }{\lim \inf }\left\vert \frac{\omega ^{\prime
}(t)}{\omega (t)^{1+\varepsilon }}\right\vert =0  \label{teeth}
\end{equation}%
Furthermore, if%
\begin{equation*}
\lim_{t\rightarrow \infty }\left\vert \frac{\omega ^{\prime }(t)}{\omega
(t)^{1+\varepsilon }}\right\vert =0
\end{equation*}%
then for all $r>0$,%
\begin{equation}
\lim_{t\rightarrow \infty }\sup \left\{ \left\vert \frac{\omega ^{\prime
}(t+s)}{\omega (t)^{1+\varepsilon }}\right\vert :s\in \mathbb{R},\left\vert
s\right\vert \leq \frac{r}{\omega (t)^{\varepsilon }}\right\} =0
\label{rho lim}
\end{equation}
\end{proposition}

A function $f:\mathbb{R}^{n}\rightarrow \mathbb{R}$ is called homothetic if
it is of the form $f(x)=\psi (h(x))$, where $\psi :\mathbb{R}\rightarrow 
\mathbb{R}$ is strictly increasing and $h:\mathbb{R}^{n}\rightarrow \mathbb{R%
}$ is positively $1$-homogeneous, i.e. $h(\alpha x)=\alpha h(x)$ for all $%
x\in \mathbb{R}^{n}$ and $a>0$, see e.g. \cite{SiBl}. For our purposes, we
do not require $\psi $ to be increasing. In the case where $f$ is the
density of a probability measure, and we assume that $\psi $ and $h$ are
measurable and that $h$ is non-negative, and we consider a random vector $X$
with density $f$, then $h(X)$ and $X/|X|$ are independent.

Let $\mathcal{M}$ be a $C^{3}$ differentiable manifold of dimension $n-1$ in 
$\mathbb{R}^{n}$ (by this, let us take as our definition that $\mathcal{M}$
is locally the graph of a $C^{3}$ function in the appropriate coordinate
system, i.e. for all $\theta \in \mathcal{M}$ there exists an open set $%
U_{\theta }\subset \mathbb{R}^{n}$, a $C^{3}$ function $f_{\theta }:\mathbb{R%
}^{n-1}\rightarrow \mathbb{R}$, and an affine isometry $I_{\theta }:\mathbb{R%
}^{n}\rightarrow \mathbb{R}^{n}$ such that $I_{\theta }\left( \theta \right)
\in U_{\theta }$ and $U_{\theta }\cap I_{\theta }\left( \mathcal{M}\right)
=U_{\theta }\cap f_{\theta }$, where we identify $f_{\theta }$ as a subset
of $\mathbb{R}^{n}$). Consider any $\theta \in \mathcal{M}$, and let $%
\mathbf{n}(\theta )$ denote one of the two unit normal vectors associated to 
$\mathcal{M}$ at $\theta $. Then there exists an affine isometry $J:\mathbb{R%
}^{n}\rightarrow \mathbb{R}^{n}$ such that $J\left( \theta \right) =0$, $%
J\left( \mathbf{n}(\theta )\right) =\left( 0,0,\ldots ,0,1\right) $, and $%
J\left( \mathcal{M}\right) $ coincides with the graph of a $C^{3}$ function $%
\psi :\mathbb{R}^{n-1}\rightarrow \mathbb{R}$ in a neighbourhood of $0\in 
\mathbb{R}^{n}$, with $\nabla \psi \left( 0\right) =0$. By the spectral
theorem applied to the Hessian $H_{\psi }\left( 0\right) $ we may choose $J$
so that $H_{\psi }\left( 0\right) $ is diagonal. The second order Taylor
expansion then takes the form%
\begin{equation*}
\psi \left( z\right) =\frac{1}{2}\sum_{i=1}^{n-1}\kappa
_{i}z_{i}^{2}+R\left( z\right)
\end{equation*}%
where $R=O\left( \left\vert z\right\vert ^{3}\right) $ as $z\rightarrow 0$.
The sequence of values $(\kappa _{i})_{1}^{n-1}$ are known as the principal
curvatures of $\mathcal{M}$ at $\theta $ with respect to $\mathbf{n}(\theta
) $, and $\left( J^{-1}e_{i}-\theta \right) _{1}^{n-1}$ is called a sequence
of principal directions. By the Gauss-Kronecker curvature of $\mathcal{M}$
at $\theta $, we mean%
\begin{equation*}
\kappa (\theta )=\det H_{\psi }\left( 0\right) =\prod_{i=1}^{n-1}\kappa _{i}
\end{equation*}%
which is independent of the choice of $J$ (but the sign of $\kappa (\theta )$
depends on the choice of $\mathbf{n}(\theta )$ when $n-1$ is odd).

These definitions, as well as Proposition \ref{pploy}, are relevant to the
following result. We label the principal curvatures $(\kappa _{i})_{2}^{n}$
instead of $(\kappa _{i})_{1}^{n-1}$.

\begin{theorem}[homotheticity]
\label{th homo level}Let $n\geq 2$ and let $h:\mathbb{R}^{n}\rightarrow 
\mathbb{R}$ be a continuous function that is not identically zero, such that
for all $x\in \mathbb{R}^{n}$ and all $\alpha >0$, $h(\alpha x)=\alpha h(x)$%
. Let $\mathcal{M}=\left\{ x\in \mathbb{R}^{n}:\left\vert h(x)\right\vert
=1\right\} $ and $\theta \in \mathcal{M}$ be such that $\mathcal{M}$ is $%
C^{3}$ is a neighbourhood of $\theta $ and such that the Gauss-Kronecker
curvature of $\mathcal{M}$ at $\theta $, denoted $\kappa (\theta )$, is
nonzero. Let $(\kappa _{i})_{2}^{n}$ be the principal curvatures of $%
\mathcal{M}$ at $\theta $ (in any order) with respect to the inward pointing
normal vector $\mathbf{n}(\theta )=-h(\theta )\left\vert \nabla h(\theta
)\right\vert ^{-1}\nabla h(\theta )$. Let $\rho :\mathbb{R}\rightarrow 
\mathbb{R}$ be any $C^{2}$ function and let $f:\mathbb{R}^{n}\rightarrow
(0,\infty )$ be defined as $f(x)=\exp \left( -\rho (h(x))\right) $. Then
there exists $C>0$ and an injective linear map $Q:\mathbb{R}%
^{n-1}\rightarrow \mathbb{R}^{n}$ such that for all $x\in \mathbb{R}^{n}$
and all $t>0$ satisfying $\rho ^{\prime }\left( th\left( \theta \right)
\right) \neq 0$,%
\begin{equation}
\frac{f\left( t\theta +T_{t}x\right) }{f(t\theta )}=\exp \left( -x_{1}-\pi
\sum_{i=2}^{n}\varepsilon _{i}x_{i}^{2}+\delta (x,t)\right)
\label{homog approx th w}
\end{equation}%
where%
\begin{eqnarray*}
\varepsilon _{i} &=&h(\theta )\frac{\rho ^{\prime }\left( th(\theta )\right) 
}{\left\vert \rho ^{\prime }\left( th(\theta )\right) \right\vert }\frac{%
\kappa _{i}}{\left\vert \kappa _{i}\right\vert } \\
T_{t}x &=&\frac{h(\theta )x_{1}}{\rho ^{\prime }(th(\theta ))}\theta +\sqrt{%
\frac{t}{\left\vert \rho ^{\prime }(th(\theta ))\right\vert }}%
Q(x_{2},x_{3},\ldots x_{n})
\end{eqnarray*}%
and $\left\vert \delta (x,t)\right\vert $ can be bounded above by%
\begin{equation*}
C\frac{\left\vert \rho ^{\prime \prime }(th(\theta )+s)\right\vert }{\rho
^{\prime }(th(\theta ))^{2}}\left[ \left\vert x_{1}\right\vert ^{2}+\left(
\sum_{i=2}^{n}\left\vert x_{i}\right\vert ^{2}\right) ^{2}+\min \left\{
I^{2},II^{2}\right\} \right] +C\min \left\{ I,II\right\}
\end{equation*}%
where%
\begin{equation*}
I=\frac{\left\vert x_{1}\right\vert ^{2}}{\left\vert t\rho ^{\prime
}(th(\theta ))\right\vert }+\sum_{i=2}^{n}\left\vert x_{i}\right\vert ^{2}%
\hspace{0.65in}II=\frac{\left\vert t\rho ^{\prime }(th(\theta ))\right\vert
^{-3/2}\left\vert x_{1}\right\vert ^{3}+\left( \sum_{i=2}^{n}\left\vert
x_{i}\right\vert ^{2}\right) ^{3/2}}{\left\vert t\rho ^{\prime }(th(\theta
))\right\vert ^{1/2}}
\end{equation*}%
and%
\begin{equation*}
\left\vert s\right\vert \leq C\left\vert \rho ^{\prime }(th(\theta
))\right\vert ^{-1}\left( \left\vert x_{1}\right\vert
+\sum_{i=2}^{n}\left\vert x_{i}\right\vert ^{2}+\min \left\{ I,II\right\}
\right)
\end{equation*}%
Lastly%
\begin{equation*}
\left\vert \mathrm{Det}(T_{t})\right\vert =\left\vert \kappa (\theta
)\right\vert ^{-1/2}\left( 2\pi t\right) ^{(n-1)/2}\left( \left\vert \nabla
h(\theta )\right\vert \rho ^{\prime }(th(\theta ))\right) ^{-(n+1)/2}
\end{equation*}%
In the special case where $h$ is a norm and $\rho $ is increasing on $%
[0,\infty )$, $\varepsilon _{i}=1$ for all $2\leq i\leq n$, and if (in
addition) $\mathcal{M}$ is $C^{3}$ everywhere with $\kappa (\theta )\neq 0$
for all $\theta \in \mathcal{M}$, then the value of $C$ may be taken
independently of $\theta $.
\end{theorem}

\begin{theorem}[mixed homogeneity]
\label{th mixed homo}Let $n,m\geq 2$ and for each $1\leq i\leq m$ let $%
p(i)\in (0,\infty )$ and let $q_{i}:\mathbb{R}^{n}\rightarrow \mathbb{R}$ be
a continuous function that is $C^{3}$ on $\mathbb{R}^{n}\backslash \{0\}$
and not identically zero, such that for all $x\in \mathbb{R}^{n}$ and all $%
\alpha \geq 0$, $q_{i}(\alpha x)=\alpha ^{p(i)}q_{i}(x)$. Assume that $%
p(1)>\max_{2\leq i\leq m}p(i)$. Define $f:\mathbb{R}^{n}\rightarrow
(0,\infty )$ as%
\begin{equation*}
f(x)=\exp \left( -\sum_{i=1}^{m}q_{i}(x)\right)
\end{equation*}%
Let $\mathcal{M}=\left\{ \theta \in \mathbb{R}^{n}:\left\vert q_{1}(\theta
)\right\vert =1\right\} $ and consider any $\theta \in \mathcal{M}$ such
that $\kappa (\theta )\neq 0$. Then there exists $(\varepsilon
_{i})_{2}^{n}\in \left\{ \pm 1\right\} ^{n-1}$ and a function $t\mapsto
T_{t} $ from $[1,\infty )$ into $GL(n)$ such that the function%
\begin{equation*}
x\mapsto \frac{f\left( t\theta +T_{t}x\right) }{f(t\theta )}
\end{equation*}%
converges uniformly on compact subsets of $\mathbb{R}^{n}$ to the function $%
\exp \left( -x_{1}-\pi \sum_{2}^{n}\varepsilon _{i}x_{i}^{2}\right) $ as $%
t\rightarrow \infty $.
\end{theorem}

\begin{theorem}[finding $T$]
\label{find scale}Let $n\geq 2$, $f:\mathbb{R}^{n}\rightarrow (0,\infty )$, $%
g=-\log f$, and let $\xi :\mathbb{R}^{n}\rightarrow \mathbb{R}$ be $C^{2}$
in a neighbourhood of $0$ with%
\begin{equation*}
\nabla \xi (0)=e_{1}\hspace{0.65in}H_{\xi }(0)=\left[ 
\begin{array}{cccc}
0 & 0 & \ldots & 0 \\ 
0 & 2\pi & \ldots & 0 \\ 
\vdots & \vdots & \ddots & \vdots \\ 
0 & 0 & \ldots & 2\pi%
\end{array}%
\right]
\end{equation*}%
Let $y\in \mathbb{R}^{n}$, and let $T$ be an invertible $n\times n$ matrix
such that for all $x\in \mathbb{R}^{n}$,%
\begin{equation}
\frac{f\left( y+Tx\right) }{f(y)}=\exp \left( -\xi (x)\right)  \label{ph}
\end{equation}%
Then $T$ can be expressed as%
\begin{equation}
T=A^{-1/2}F^{T}  \label{gives T}
\end{equation}%
where $A$ is the positive definite $n\times n$ matrix defined by 
\begin{equation*}
A_{ij}=\frac{1}{2\pi }\partial _{ij}g(y)+\partial _{i}g(y)\partial _{j}g(y)
\end{equation*}%
$A^{-1/2}$ is the principal square root of $A^{-1}$, and $F\in O(n)$ such
that%
\begin{equation}
FA^{-1/2}\nabla g(y)=e_{1}  \label{ro1}
\end{equation}%
This also implies, by (\ref{gives T}) and (\ref{ro1}) that%
\begin{equation}
Te_{1}=A^{-1}\nabla g(y)  \label{aff grad}
\end{equation}%
Furthermore, in the special case where $\xi (x)=x_{1}+\pi
\sum_{i=2}^{n}x_{i}^{2}$, for any $G\in O(n)$ such that%
\begin{equation}
GA^{-1/2}\nabla g(y)=e_{1}  \label{ro2}
\end{equation}%
it follows that for all $x\in \mathbb{R}^{n}$,%
\begin{equation*}
\frac{f\left( y+A^{-1/2}G^{T}x\right) }{f(y)}=\exp \left( -x_{1}-\pi
\sum_{i=2}^{n}x_{i}^{2}\right)
\end{equation*}
\end{theorem}

\begin{theorem}[homogeneity of $\Lambda _{n}$]
\label{self sim}Let $n\geq 2$. Then for all $y\in \mathbb{R}^{n}$ there
exists an invertible $n\times n$ matrix $T$ with the following property. For
all $x\in \mathbb{R}^{n}$,%
\begin{equation*}
\frac{\Lambda _{n}\left( y+Tx\right) }{\Lambda _{n}(y)}=\Lambda _{n}(x)
\end{equation*}
\end{theorem}

\section{Proofs}

\begin{proof}[Proof of Theorem \protect\ref{extremal proc}]
The proof is based on the well known concept of Poissonization: the
empirical distribution of an i.i.d. random sample of random size is a
Poisson point process, as long as the size of the sample follows a Poisson
distribution and is independent of the sample points themselves. Since the
underlying probability space $\Omega $ is not arbitrary, we may assume that
it is rich enough so that we may introduce new independent random objects as
necessary. Let us assume, as we may, that the random vectors $%
(X_{i})_{1}^{\infty }$ are generated as follows: Note that the set $\left\{
(x,t)\in \mathbb{R}^{n}\times \lbrack 0,\infty ):t\leq f(x)\right\} $ has
Lebesgue measure $1$, and let $(X_{i}^{\ast })_{1}^{\infty }$ denote an
i.i.d. sequence of random vectors uniformly distributed in this set. Then
set $X_{i}=PX_{i}^{\ast }$, where $P:\mathbb{R}^{n}\times \lbrack 0,\infty
)\rightarrow \mathbb{R}^{n}$ is the natural projection defined by $P(x,t)=x$%
. This construction is useful because it will allow us to decompose a
certain Poisson process into the sum of two (particular) independent Poisson
processes. For each $i$ the random vector $T^{-1}(X_{i}-y)$ has a
distribution with density $\det (T)f\left( y+Tx\right) $. Defining $Q:%
\mathbb{R}^{n}\times \mathbb{R}\rightarrow \mathbb{R}^{n}\times \mathbb{R}$
as%
\begin{equation*}
Q(x,t)=\left( T^{-1}(x-y),N\det (T)t\right)
\end{equation*}%
we see that each $QX_{i}^{\ast }$ is uniformly distributed in%
\begin{equation*}
\left\{ (x,t)\in \mathbb{R}^{n}\times \lbrack 0,\infty ):0\leq t\leq N\det
(T)f\left( y+Tx\right) \right\}
\end{equation*}%
and that $PQX_{i}^{\ast }=T^{-1}(X_{i}-y)$. Now define $W_{1}$ as%
\begin{equation*}
\left\{ (x,t)\in B^{\circ }(0,R)\times \lbrack 0,\infty ):0\leq t\leq \min
\left\{ N\det (T)f\left( y+Tx\right) ,\Lambda \left( x\right) \right\}
\right\}
\end{equation*}%
and%
\begin{eqnarray*}
W_{2} &=&\left\{ (x,t)\in B^{\circ }(0,R)\times \lbrack 0,\infty ):\Lambda
\left( x\right) <t\leq N\det (T)f\left( y+Tx\right) \right\} \\
W_{3} &=&\left\{ (x,t)\in B^{\circ }(0,R)\times \lbrack 0,\infty ):N\det
(T)f\left( y+Tx\right) <t\leq \Lambda \left( x\right) \right\}
\end{eqnarray*}%
Consider a Poisson distributed random variable $M\sim Pois(N)$ independent
of the sequence $(X_{i}^{\ast })_{1}^{\infty }$. The random measures on $%
\mathbb{R}^{n+1}$ defined by%
\begin{eqnarray*}
\nu _{1} &=&\sum_{i=1}^{M}\delta \left( QX_{i}^{\ast }\right) \mathbb{I}%
\left\{ QX_{i}^{\ast }\in W_{1}\right\} \\
\nu _{2} &=&\sum_{i=1}^{M}\delta \left( QX_{i}^{\ast }\right) \mathbb{I}%
\left\{ QX_{i}^{\ast }\in W_{2}\right\}
\end{eqnarray*}%
are then Poisson point processes with a constant (density of) intensity of $%
1 $ in $W_{1}$ and $W_{2}$ respectively, where $\mathbb{I}(A)$ denotes the
indicator function of an event $A$. Now let $\nu _{3}$ be a Poisson point
process of intensity $1$ in $W_{3}$ independent of all other random
quantities introduced so far. Then the pushforward measure $\nu =P(\nu
_{1}+\nu _{3})$ is a Poisson process in $B^{\circ }(0,R)$ with (density of)
intensity $\exp \left( -x_{1}-\pi \sum_{2}^{n}x_{i}^{2}\right) $. By
assumption, $\mathrm{vol}_{n+1}\left( W_{2}\cup W_{3}\right) <\varepsilon $,
and therefore by definition of a Poisson process, with probability at least $%
\exp (-\varepsilon )\geq 1-\varepsilon $, $\nu _{2}=\nu _{3}=0$. Let this
event be denoted $E_{1}$. By concentration of the Poisson distribution, see
e.g. \cite{JLR}, with probability at least $1-2N^{-1/4}$, $\left\vert
N-M\right\vert \leq \sqrt{N\log N}$ (this event to be denoted $E_{2}$). With
probability at least $1-3\sqrt{N\log N}\mu \left( y+TB^{\circ }(0,R)\right) $%
, for all $i$ with $\left\vert N-i\right\vert \leq \sqrt{N\log N}$, $%
X_{i}\notin y+TB^{\circ }(0,R)$ (this event $E_{3}$). Note that%
\begin{eqnarray*}
&&\mu \left( y+TB^{\circ }(0,R)\right) \\
&=&\int_{\left\vert x\right\vert <R}\det (T)f\left( y+Tx\right) dx\leq
N^{-1}\left( 1+\int_{\left\vert x\right\vert <R}\exp \left( -x_{1}-\pi
\sum_{2}^{n}x_{i}^{2}\right) dx\right)
\end{eqnarray*}%
Now suppose that these three events occur. Then for any Borel set $%
E\subseteq B^{\circ }(0,R)$, 
\begin{eqnarray*}
\nu _{N}(E) &=&\left\vert \left\{ i\leq N:QX_{i}^{\ast }\in E\times \lbrack
0,\infty )\right\} \right\vert =\nu _{1}(E\times \lbrack 0,\infty ))+\nu
_{2}(E\times \lbrack 0,\infty )) \\
&=&\nu _{1}(E\times \lbrack 0,\infty ))+\nu _{3}(E\times \lbrack 0,\infty
))=\nu (E)
\end{eqnarray*}%
and the proposition is proved.
\end{proof}

\textbf{Note:} The function $\exp \left( -x_{1}-\pi
\sum_{2}^{n}x_{i}^{2}\right) $ plays only a superficial role in Proposition %
\ref{extremal proc} and can be replaced with a more general function. Also,
the probability bound can obviously be improved if we replace $\nu
_{N}(E)=\nu (E)$ with $\left\vert \nu _{N}(E)-\nu (E)\right\vert \leq k$.

The following lemma is well known, but we include its short proof for
completeness. See for example Corollary 8.17 in \cite{Gard} for essentially
the same statement expressed in a different way.

\begin{lemma}
\label{volofsec}Let $H$ be a real symmetric positive definite matrix and $%
\mathcal{E}=H^{-1}(B_{2}^{n})$. Then for all $\theta \in S^{n-1}$,%
\begin{equation*}
\mathrm{vol}_{n-1}\left( \mathcal{E}\cap \theta ^{\bot }\right) =\frac{%
\mathrm{vol}_{n-1}(B_{2}^{n-1})}{\left\vert H^{-1}\theta \right\vert \mathrm{%
Det}(H)}
\end{equation*}
\end{lemma}

\begin{proof}
Let $\phi _{n}$ denote the parallel section function of $B_{2}^{n}$ (in any
direction),%
\begin{equation*}
\phi _{n}(t)=\mathrm{vol}_{n-1}\left\{ x\in B_{2}^{n}:x_{1}=t\right\}
\end{equation*}%
and let $\psi $ denote the parallel section function of $\mathcal{E}$ in the
direction of $\theta $,%
\begin{equation*}
\psi (t)=\mathrm{vol}_{n-1}\left( \mathcal{E}\cap \left( \theta ^{\bot
}+t\theta \right) \right)
\end{equation*}%
Since a non-degenerate ellipsoid is a Euclidean ball when viewed in the
correct coordinate system, and all parallel section functions of the
Euclidean ball are the same, $\psi (t)=a\phi _{n}(bt)$ for some $a,b>0$.
Recall that the dual Minkowski functional of $\mathcal{E}$ is defined as $%
\left\Vert y\right\Vert _{\mathcal{E}^{\circ }}=\sup \left\{ \left\langle
x,y\right\rangle :x\in \mathcal{E}\right\} $. Since $\psi (0)=\mathrm{vol}%
_{n-1}\left( \mathcal{E}\cap \{\theta \}^{\bot }\right) $ and the support of 
$\psi $ is the interval $\left[ -\left\Vert \theta \right\Vert _{\mathcal{E}%
^{\circ }},\left\Vert \theta \right\Vert _{\mathcal{E}^{\circ }}\right] $
whereas the support of $\phi _{n}$ is $\left[ -1,1\right] $, it follows that%
\begin{equation*}
a=\mathrm{vol}_{n-1}\left( \mathcal{E}\cap \{\theta \}^{\bot }\right) /%
\mathrm{vol}_{n-1}(B_{2}^{n-1})
\end{equation*}%
and $b=\left\Vert \theta \right\Vert _{\mathcal{E}^{\circ }}^{-1}$. Therefore%
\begin{equation*}
\mathrm{vol}_{n}\left( \mathcal{E}\right) =\int_{-\left\Vert \theta
\right\Vert _{\mathcal{E}^{\circ }}}^{\left\Vert \theta \right\Vert _{%
\mathcal{E}^{\circ }}}\psi (t)dt=\left\Vert \theta \right\Vert _{\mathcal{E}%
^{\circ }}\times \frac{\mathrm{vol}_{n-1}\left( \mathcal{E}\cap \theta
^{\bot }\right) }{\mathrm{vol}_{n-1}(B_{2}^{n-1})}\times \mathrm{vol}%
_{n}(B_{2}^{n})
\end{equation*}%
The result now follows because%
\begin{equation*}
\left\Vert \theta \right\Vert _{\mathcal{E}^{\circ }}=\sup \left\{
\left\langle \theta ,H^{-1}x\right\rangle :x\in B_{2}^{n}\right\} =\sup
\left\{ \left\langle H^{-1}\theta ,x\right\rangle :x\in B_{2}^{n}\right\}
=\left\vert H^{-1}\theta \right\vert
\end{equation*}
\end{proof}

\begin{proof}[Proof of Theorem \protect\ref{th log concave product}]
Consider the following inner product and corresponding norm defined on $%
\mathbb{R}^{n}$,%
\begin{equation*}
\left\langle u,v\right\rangle _{\sharp }=\frac{1}{2\pi }%
\sum_{i=1}^{n}g_{i}^{\prime \prime }(y_{i})u_{i}v_{i}\text{\hspace{0.65in}}%
||u||_{\sharp }=\left( \frac{1}{2\pi }\sum_{i=1}^{n}g_{i}^{\prime \prime
}(y_{i})u_{i}^{2}\right) ^{1/2}
\end{equation*}%
Since any two real Hilbert spaces of dimension $n-1$ are linearly isometric,
there exists a linear map $Q:\mathbb{R}^{n-1}\rightarrow \nabla f(y)^{\bot }$
such that $||Qz||_{\sharp }=|z|$ for all $z\in \mathbb{R}^{n-1}$. Define $T:%
\mathbb{R}^{n}\rightarrow \mathbb{R}^{n}$ by%
\begin{equation*}
Tx=\frac{\nabla g(y)}{|\nabla g(y)|^{2}}x_{1}+Q(x_{2},x_{3},\ldots x_{n})
\end{equation*}%
Now fix any $x\in \mathbb{R}^{n}$ and define%
\begin{equation*}
\alpha =f(y)^{-1}\hspace{0.65in}u=\frac{\nabla g(y)}{|\nabla g(y)|^{2}}%
\hspace{0.65in}v=Q(x_{2},x_{3},\ldots x_{n})
\end{equation*}%
Introduce the parameter $s\in \mathbb{R}$ and define%
\begin{eqnarray*}
\psi (s) &=&-\log \alpha f(y+sTx)=\log
f(y)+\sum_{i=1}^{n}g_{i}(y_{i}+su_{i}x_{1}+sv_{i}) \\
\eta (s) &=&-\log \Lambda _{n}(sx)=sx_{1}+\pi s^{2}\sum_{i=2}^{n}x_{i}^{2}
\end{eqnarray*}%
Our goal is to estimate%
\begin{equation*}
\left\vert \log \frac{f(y+Tx)}{f(y)\Lambda _{n}(x)}\right\vert =\left\vert
\psi (1)-\eta (1)\right\vert
\end{equation*}%
which by Taylor's theorem can be bounded above by%
\begin{equation*}
\left\vert \psi (0)-\eta (0)\right\vert +\left\vert \psi ^{\prime }(0)-\eta
^{\prime }(0)\right\vert +\frac{1}{2}\left\vert \psi ^{\prime \prime
}(0)-\eta ^{\prime \prime }(0)\right\vert +\frac{1}{6}\left\vert \psi
^{\prime \prime \prime }(\xi )-\eta ^{\prime \prime \prime }(\xi )\right\vert
\end{equation*}%
for some $\xi \in (0,1)$. Since $v\in Q(\mathbb{R}^{n-1})=\nabla g(y)^{\bot
} $,%
\begin{eqnarray*}
\psi ^{\prime }(s) &=&\sum_{i=1}^{n}g_{i}^{\prime
}(y_{i}+su_{i}x_{1}+sv_{i})(u_{i}x_{1}+v_{i}) \\
\psi ^{\prime }(0) &=&x_{1}\left\langle \nabla g(y),u\right\rangle
+\left\langle \nabla g(y),v\right\rangle =x_{1}
\end{eqnarray*}%
and%
\begin{eqnarray*}
\psi ^{\prime \prime }(s) &=&\sum_{i=1}^{n}g_{i}^{\prime \prime
}(y_{i}+su_{i}x_{1}+sv_{i})(u_{i}x_{1}+v_{i})^{2} \\
\psi ^{\prime \prime }(0) &=&x_{1}^{2}\sum_{i=1}^{n}g_{i}^{\prime \prime
}(y_{i})\frac{g_{i}^{\prime }(y_{i})^{2}}{|\nabla g(y)|^{4}}%
+2x_{1}\sum_{i=1}^{n}g_{i}^{\prime \prime }(y_{i})\frac{g_{i}^{\prime
}(y_{i})}{|\nabla g(y)|^{2}}v_{i}+\sum_{i=1}^{n}g_{i}^{\prime \prime
}(y_{i})v_{i}^{2} \\
&=&x_{1}^{2}\sum_{i=1}^{n}\frac{g_{i}^{\prime \prime }(y_{i})}{g_{i}^{\prime
}(y_{i})^{2}}\frac{g_{i}^{\prime }(y_{i})^{4}}{|\nabla g(y)|^{4}}%
+2x_{1}\sum_{i=1}^{n}\frac{g_{i}^{\prime \prime }(y_{i})^{1/2}}{%
g_{i}^{\prime }(y_{i})}\frac{g_{i}^{\prime }(y_{i})^{2}}{|\nabla g(y)|^{2}}%
g_{i}^{\prime \prime }(y_{i})^{1/2}v_{i} \\
&&+\sum_{i=1}^{n}g_{i}^{\prime \prime }(y_{i})v_{i}^{2}
\end{eqnarray*}%
where the last term may be written as%
\begin{equation*}
\sum_{i=1}^{n}g_{i}^{\prime \prime }(y_{i})v_{i}^{2}=2\pi ||v||_{\sharp
}^{2}=2\pi \sum_{i=2}^{n}x_{i}^{2}=\eta ^{\prime \prime }(0)
\end{equation*}%
and therefore%
\begin{eqnarray*}
&&\frac{1}{2}\left\vert \psi ^{\prime \prime }(0)-\eta ^{\prime \prime
}(0)\right\vert \\
&\leq &\frac{1}{2}\left\vert x_{1}\right\vert ^{2}\sum_{i=1}^{n}\frac{%
g_{i}^{\prime \prime }(y_{i})}{g_{i}^{\prime }(y_{i})^{2}}\frac{%
g_{i}^{\prime }(y_{i})^{4}}{|\nabla g(y)|^{4}}+\sqrt{2\pi }\left\vert
x_{1}\right\vert \left( \sum_{i=2}^{n}x_{i}^{2}\right) ^{1/2}\left(
\sum_{i=1}^{n}\frac{g_{i}^{\prime \prime }(y_{i})}{g_{i}^{\prime }(y_{i})^{2}%
}\frac{g_{i}^{\prime }(y_{i})^{4}}{|\nabla g(y)|^{4}}\right) ^{1/2}
\end{eqnarray*}%
Similarly,%
\begin{eqnarray*}
\psi ^{\prime \prime \prime }(s) &=&\sum_{i=1}^{n}g_{i}^{\prime \prime
\prime }(y_{i}+su_{i}x_{1}+sv_{i})(u_{i}x_{1}+v_{i})^{3} \\
&=&x_{1}^{3}\sum_{i=1}^{n}\frac{g_{i}^{\prime \prime \prime
}(y_{i}+su_{i}x_{1}+sv_{i})}{g_{i}^{\prime \prime }(y_{i})^{3/2}}\frac{%
g_{i}^{\prime }(y_{i})^{6}}{|\nabla g(y)|^{6}}\frac{g_{i}^{\prime \prime
}(y_{i})^{3/2}}{g_{i}^{\prime }(y_{i})^{3}} \\
&&+3x_{1}^{2}\sum_{i=1}^{n}\frac{g_{i}^{\prime \prime \prime
}(y_{i}+su_{i}x_{1}+sv_{i})}{g_{i}^{\prime \prime }(y_{i})^{3/2}}%
g_{i}^{\prime \prime }(y_{i})^{1/2}v_{i}\frac{g_{i}^{\prime }(y_{i})^{4}}{%
|\nabla g(y)|^{4}}\frac{g_{i}^{\prime \prime }(y_{i})}{g_{i}^{\prime
}(y_{i})^{2}} \\
&&+3x_{1}\sum_{i=1}^{n}\frac{g_{i}^{\prime \prime \prime
}(y_{i}+su_{i}x_{1}+sv_{i})}{g_{i}^{\prime \prime }(y_{i})^{3/2}}%
g_{i}^{\prime \prime }(y_{i})v_{i}^{2}\frac{g_{i}^{\prime }(y_{i})^{2}}{%
|\nabla g(y)|^{2}}\frac{g_{i}^{\prime \prime }(y_{i})^{1/2}}{g_{i}^{\prime
}(y_{i})} \\
&&+\sum_{i=1}^{n}\frac{g_{i}^{\prime \prime \prime
}(y_{i}+su_{i}x_{1}+sv_{i})}{g_{i}^{\prime \prime }(y_{i})^{3/2}}%
(g_{i}^{\prime \prime }(y_{i})^{1/2}v_{i})^{3}
\end{eqnarray*}%
which can be bounded above as before. Lastly, using Lemma \ref{volofsec}
with $\theta =\left\vert \nabla g(y)\right\vert ^{-1}\nabla g(y)$ and $%
H=(2\pi )^{-1/2}H_{g}(y)^{1/2}$, and observing that $Q(B_{2}^{n-1})=\mathcal{%
E}\cap \nabla g(y)^{\bot }$, 
\begin{equation*}
\left\vert \mathrm{Det}(T)\right\vert =\left\vert \nabla g(y)\right\vert
^{-1}\times \frac{\mathrm{vol}_{n-1}(Q(B_{2}^{n-1}))}{\mathrm{vol}%
_{n-1}(B_{2}^{n-1})}=\frac{1}{\left\vert H^{-1}\nabla g(y)\right\vert \cdot 
\mathrm{Det}(H)}
\end{equation*}
\end{proof}

\begin{proof}[Proof of Proposition \protect\ref{pploy}]
Let $\varepsilon $ and $\omega $ be as in the statement of the proposition,
and suppose that the $\lim \inf $ in (\ref{teeth}) is nonzero. Then there
exists $t_{0},c>0$ such that for all $t\geq t_{0}$, $\omega (t)>0$ and $%
|\omega ^{\prime }(t)|>c\omega (t)^{1+\varepsilon }$. In particular, $\omega
^{\prime }(t)\neq 0$ and $\omega $ is either strictly increasing on $%
[t_{0},\infty )$, or strictly decreasing. Hence the sign of $\omega ^{\prime
}$ is constant on $[t_{0},\infty )$. It also follows that $\omega $ is
injective on $[t_{0},\infty )$ and therefore satisfies an autonomous
differential equation $\omega ^{\prime }(t)=\Theta (\omega (t))$, where $%
\Theta (s)=\omega ^{\prime }(\omega ^{-1}(s))$ and $\left\vert \Theta
(s)\right\vert >cs^{1+\varepsilon }$ for all $s\in \left\{ \omega (t):t\in
\lbrack t_{0},\infty )\right\} $. In what follows it may help the reader to
consider two cases separately: $\omega ^{\prime }(t_{0})<0$ and $\omega
^{\prime }(t_{0})>0$. Consider the interval%
\begin{equation*}
I=\left\{ t\in \lbrack t_{0},\infty ):\omega (t_{0})^{-\varepsilon }-\frac{%
c\varepsilon \omega ^{\prime }(t_{0})\left( t-t_{0}\right) }{\left\vert
\omega ^{\prime }(t_{0})\right\vert }>0\right\}
\end{equation*}%
We now compare $\omega $ to the function%
\begin{equation*}
u(t)=\left( \omega (t_{0})^{-\varepsilon }-\frac{c\varepsilon \omega
^{\prime }(t_{0})\left( t-t_{0}\right) }{\left\vert \omega ^{\prime
}(t_{0})\right\vert }\right) ^{-1/\varepsilon }\text{, }t\in I
\end{equation*}%
which is the solution to the initial value problem%
\begin{equation*}
u(t_{0})=\omega (t_{0})\hspace{0.65in}u^{\prime }(t)=c\frac{\omega ^{\prime
}(t_{0})u(t)^{1+\varepsilon }}{\left\vert \omega ^{\prime
}(t_{0})\right\vert }\text{, }t\in I
\end{equation*}%
The functions $\omega $ and $u$ therefore satisfy%
\begin{equation*}
\int_{\omega (t_{0})}^{\omega (t)}\frac{1}{\Theta (s)}ds=t-t_{0}\hspace{%
0.65in}\int_{u(t_{0})}^{u(t)}\frac{\omega ^{\prime }(t_{0})}{c\left\vert
\omega ^{\prime }(t_{0})\right\vert s^{1+\varepsilon }}ds=t-t_{0}
\end{equation*}%
(just differentiate both sides), where we use a standard convention that $%
\int_{a}^{b}=-\int_{b}^{a}$. This implies that on $I$, $\omega ^{\prime
}(t_{0})\left( \omega (t)-u(t)\right) \geq 0$, for if $\omega (a)=u(b)$ then 
$a\leq b$. This fact, which is a standard monotonicity property of
autonomous ODEs, means that $\omega $ outruns $u$. If $\omega ^{\prime
}(t_{0})<0$ then $\omega (t)\leq u(t)$ on $[t_{0},\infty )$, which
contradicts our assumption that $t^{1/\varepsilon }\omega (t)\rightarrow
\infty $. If $\omega ^{\prime }(t_{0})>0$ then $\omega (t)\geq u(t)$ on $I$
which is also a contradiction since in this case $u$ explodes to infinity in
finite time. We now move on to the second part of the proposition. Since $%
t^{1/\varepsilon }\omega (t)\rightarrow \infty $ as $t\rightarrow \infty $,
there exists $t_{0}>0$ such that $\omega (t)>0$ for all $t\geq t_{0}$. Fix
any $r,\delta >0$ and consider the initial value problem%
\begin{equation*}
\psi (0)=t_{0}\hspace{0.65in}\psi ^{\prime }(x)=\omega \left( \psi
(x)\right) ^{-\varepsilon }\text{, }x>0
\end{equation*}%
A unique solution to this differential equation exists on $[0,\infty )$, by
the basic theory of autonomous ODEs, and%
\begin{equation*}
\left\vert \frac{d}{dx}\log \psi ^{\prime }(x)\right\vert =\frac{\varepsilon
\left\vert \omega ^{\prime }\left( \psi (x)\right) \right\vert }{\omega
\left( \psi (x)\right) ^{1+\varepsilon }}
\end{equation*}%
It also follows that $\psi $ is strictly increasing and $\lim_{t\rightarrow
\infty }\psi (t)=\infty $ (if not then $\lim_{t\rightarrow \infty }\psi
(t)=a<\infty $ and $\lim_{t\rightarrow \infty }\psi ^{\prime }(t)=\omega
\left( a\right) ^{-\varepsilon }>0$ which in turn implies that $%
\lim_{t\rightarrow \infty }\psi (t)=\infty $). By assumption, there exists $%
x_{0}>t_{0}$ such that for all $x>x_{0}$,%
\begin{equation}
\frac{\left\vert \omega ^{\prime }\left( x\right) \right\vert }{\omega
\left( x\right) ^{1+\varepsilon }}<\min \left\{ 1,\varepsilon ^{-1}\right\}
\times \min \left\{ 1,\exp \left( -\frac{1+\varepsilon }{50\varepsilon }%
\right) \right\} \times \min \left\{ \frac{1}{100r},\frac{\delta }{100}%
\right\}  \label{agbp}
\end{equation}%
and there exists $t_{2}>x_{0}$ such that for all $t>t_{2}$, $\psi \left(
\psi ^{-1}(t)-2r\right) >x_{0}$. Now consider any $t>t_{2}$ and $s\in 
\mathbb{R}$ such that $\left\vert s\right\vert \leq r/\omega
(t)^{\varepsilon }$. It follows that for all $x\geq \psi ^{-1}(t)-2r$,%
\begin{equation*}
\left\vert \frac{d}{dx}\log \psi ^{\prime }(x)\right\vert <\frac{1}{100r}
\end{equation*}%
which implies that $\left\vert \log \psi ^{\prime }(x)-\log \psi ^{\prime
}(\psi ^{-1}(t))\right\vert <\left\vert x-\psi ^{-1}(t)\right\vert /(100r)$,
which can be rewritten as%
\begin{equation}
\exp \left( \frac{-1}{100r}\left\vert x-\psi ^{-1}(t)\right\vert \right)
\psi ^{\prime }(\psi ^{-1}(t))<\psi ^{\prime }(x)<\exp \left( \frac{1}{100r}%
\left\vert x-\psi ^{-1}(t)\right\vert \right) \psi ^{\prime }(\psi ^{-1}(t))
\label{ahhhhp}
\end{equation}%
Since $\exp (1/50)\approx 1.02$, 
\begin{eqnarray*}
\psi \left( \psi ^{-1}(t)-2r\right) &\leq &\psi (\psi ^{-1}(t))-1.8r\psi
^{\prime }(\psi ^{-1}(t))=t-1.8\frac{r}{\omega (t)^{\varepsilon }} \\
\psi \left( \psi ^{-1}(t)+2r\right) &\geq &\psi (\psi ^{-1}(t))+1.8r\psi
^{\prime }(\psi ^{-1}(t))=t+1.8\frac{r}{\omega (t)^{\varepsilon }}
\end{eqnarray*}%
which implies%
\begin{eqnarray*}
\psi \left( \psi ^{-1}(t)-2r\right) &\leq &t+s\leq \psi \left( \psi
^{-1}(t)+2r\right) \\
\psi ^{-1}(t)-2r &\leq &\psi ^{-1}(t+s)\leq \psi ^{-1}(t)+2r
\end{eqnarray*}%
Lastly, by (\ref{agbp}) and (\ref{ahhhhp}),%
\begin{equation*}
\frac{\left\vert \omega ^{\prime }\left( t+s\right) \right\vert }{\omega
\left( t\right) ^{1+\varepsilon }}=\frac{\left\vert \omega ^{\prime }\left(
t+s\right) \right\vert }{\omega \left( t+s\right) ^{1+\varepsilon }}\left( 
\frac{\psi ^{\prime }\left( \psi ^{-1}(t+s)\right) }{\psi ^{\prime }\left(
\psi ^{-1}(t)\right) }\right) ^{-(1+\varepsilon )/\varepsilon }\leq \frac{%
\delta }{100}
\end{equation*}%
Since this holds for all such $s$ and $\delta $, (\ref{rho lim}) follows.
\end{proof}

\begin{proof}[Proof of Theorem \protect\ref{th homo level}]
The first part of the proof is to establish (\ref{local taylor}) with error
term $\left\vert R(x)\right\vert \leq C\min \left\{ \left\vert x\right\vert
^{2},\left\vert x\right\vert ^{3}\right\} $. It may happen that such an
expansion exists without the assumption that $\mathcal{M}$ is $C^{3}$ near $%
\theta $, and in that case the conclusion of the theorem still holds. The
existence of a Taylor expansion of $h$ around $\theta $ follows immediately
from the $C^{3}$ condition. Our goal is to show that it takes a particular
form, as in (\ref{local taylor}), i.e. essentially to compute the Hessian $%
H_{h}\left( \theta \right) $. It is our preference to give full detail here.
This first part of the proof is more of a technicality, and the reader may
choose to skip directly to (\ref{local taylor}) onwards to see how the
approximation (\ref{homog approx th w}) comes about. By homogeneity and
continuity of $h$, $\mathcal{M}$ is homeomorphic to a subset $V\subseteq
S^{n-1}$ (open in the subspace topology of $S^{n-1}$). This homeomorphism $%
\varphi :\mathcal{M}\rightarrow V$ and its inverse are given by%
\begin{equation*}
\varphi (x)=\left\vert x\right\vert ^{-1}x\hspace{0.65in}\varphi
^{-1}(x)=\left\vert h(x)\right\vert ^{-1}x
\end{equation*}%
Since $\mathcal{M}$ is assumed to be $C^{3}$ in an $(n-1)$-dimensional
neighbourhood of $\theta \in \mathcal{M}$, the radial function $%
rad:V\rightarrow (0,\infty )$ defined by $rad(x)=\left\vert \varphi
^{-1}(x)\right\vert =\left\vert h(x)\right\vert ^{-1}$ is $C^{3}$ in an $%
(n-1)$-dimensional neighbourhood of $\left\vert \theta \right\vert
^{-1}\theta \in V$, and therefore $h(x)=\left\vert x\right\vert h\left(
\left\vert x\right\vert ^{-1}x\right) =\left\vert x\right\vert \left(
rad\left( \left\vert x\right\vert ^{-1}x\right) \right) ^{-1}$ is $C^{3}$ in
an $n$-dimensional neighbourhood of $\theta \in \mathbb{R}^{n}$. The
homogeneity condition and the fact that $h\left( \theta \right) \neq 0$
ensure that $\nabla h\left( \theta \right) \neq 0$. Consider the
inward-pointing normal vector $\mathbf{n}(\theta )=-h(\theta )\nabla
h(\theta )/\left\vert \nabla h(\theta )\right\vert $. Let $(v_{i})_{2}^{n}$
denote a sequence of principal directions of $\mathcal{M}$ at $\theta $ and $%
(\kappa _{i})_{2}^{n}$ the corresponding principal curvatures (defined with
respect to $\mathbf{n}(\theta )$). Let $W:\mathbb{R}^{n-1}\rightarrow \nabla
h(\theta )^{\perp }$ be the linear isometry (with respect to the standard
Euclidean metric on both spaces) such that $We_{i}=v_{i+1}$ $(1\leq i\leq
n-1)$. It follows that there exists a connected neighbourhood $U$ of $0\in 
\mathbb{R}^{n-1}$ and a continuous function $\psi :U\rightarrow \mathbb{R}$
of the form%
\begin{equation*}
\psi (z)=\frac{1}{2}\sum_{i=1}^{n-1}\kappa _{i+1}z_{i}^{2}+R_{1}(z)
\end{equation*}%
where $R_{1}(z)=O\left( \left\vert z\right\vert ^{3}\right) $ as $%
z\rightarrow 0$, such that for all $z\in U$, $\theta +W(z)+\psi (z)\mathbf{n}%
(\theta )\in \mathcal{M}$. We may assume (after possibly choosing a smaller $%
U$), that $h(\theta +W(z))/h(\theta )>0$ for all $z\in U$. Let%
\begin{eqnarray*}
\eta (z) &=&1-\frac{h(\theta )}{h\left( \theta +W(z)\right) } \\
a &=&\theta +W(z)+\psi (z)\mathbf{n}(\theta ) \\
b &=&\left( 1-\eta (z)\right) (\theta +W(z))
\end{eqnarray*}%
By the mean value theorem there exists $\lambda \in (0,1)$ such that,
setting $\xi =\lambda a+(1-\lambda )b$, 
\begin{equation}
h(b)=h\left( a\right) +\left\langle b-a,\nabla h\left( \xi \right)
\right\rangle  \label{long}
\end{equation}%
By homogeneity of $h$ and by definition of $\eta (\cdot )$, $h(b)=h(\theta )$%
. Since $a\in \mathcal{M}$, $h(a)\in \left\{ \pm 1\right\} $. However since $%
\psi $ and $h$ are continuous and $U$ is connected, and $\psi (0)=0$, it
follows that $h(a)=h(\theta )$. Therefore (\ref{long}) can be rewritten as%
\begin{equation*}
\eta (z)\left\langle \theta +W(z),\nabla h\left( \xi \right) \right\rangle
=\psi (z)h(\theta )\left\langle \frac{\nabla h(\theta )}{\left\vert \nabla
h(\theta )\right\vert },\nabla h\left( \xi \right) \right\rangle
\end{equation*}%
As $z\rightarrow 0$, $W(z)\rightarrow 0$ and $\nabla h\left( \xi \right)
\rightarrow \nabla h\left( \theta \right) $ (since $h$ is $C^{1}$ at $\theta 
$). By definition of $W$, $\left\langle W(z),\nabla h(\theta )\right\rangle
=0$, and by homogeneity $\left\langle \theta ,\nabla h\left( \theta \right)
\right\rangle =h(\theta )$. Therefore as $z\rightarrow 0$,%
\begin{equation*}
\eta (z)=(1+o(1))\psi (z)\left\vert \nabla h(\theta )\right\vert
\end{equation*}%
and%
\begin{eqnarray}
h\left( \theta +W(z)\right) &=&\frac{h(\theta )}{1-\eta (z)}=h(\theta
)\left( 1+(1+o(1))\eta (z)\right)  \notag \\
&=&h(\theta )+\frac{1}{2}h(\theta )\left\vert \nabla h(\theta )\right\vert
\sum_{i=1}^{n-1}\kappa _{i+1}z_{i}^{2}+R_{2}(z)  \label{some approxe}
\end{eqnarray}%
where $R_{2}(z)=o\left( \left\vert z\right\vert ^{2}\right) $ as $%
z\rightarrow 0$ (the error term will be improved to $O\left( \left\vert
z\right\vert ^{3}\right) $ in a moment). Consider the $(n-1)\times (n-1)$
diagonal matrix $E$ defined by%
\begin{equation*}
E_{i,i}=\sqrt{\frac{2\pi }{\left\vert \nabla h(\theta )\right\vert \cdot
\left\vert \kappa _{i+1}\right\vert }}
\end{equation*}%
and set $Q=WE$, in which case (\ref{some approxe}) transforms to%
\begin{equation}
h\left( \theta +Qz\right) =h(\theta )+\pi \sum_{i=1}^{n-1}\widehat{%
\varepsilon }_{i+1}z_{i}^{2}+R_{3}(z)  \label{pre Taylor}
\end{equation}%
where $\widehat{\varepsilon }_{i}=h(\theta )\kappa _{i}/\left\vert \kappa
_{i}\right\vert $. This Taylor approximation in $(n-1)$ variables extends to
an approximation in $n$ variables by homogeneity, as follows. Recycling the
variable $z$, set%
\begin{equation*}
z=\theta +x_{1}h(\theta )\theta +Q(x_{2},x_{3},\ldots x_{n})
\end{equation*}%
Assuming without loss of generality that $\left\vert x_{1}\right\vert <1/2$,
set $s=\left( 1+x_{1}h(\theta )\right) ^{-1}$. By homogeneity of $h$, this
definition of $z$, and using (\ref{pre Taylor}),%
\begin{eqnarray*}
h(z) &=&\left( 1+x_{1}h(\theta )\right) h\left( \left( 1+x_{1}h(\theta
)\right) ^{-1}z\right) \\
&=&\left( 1+x_{1}h(\theta )\right) h\left( \theta +Q(sx_{2},sx_{3},\ldots
sx_{n})\right) \\
&=&\left( 1+x_{1}h(\theta )\right) \left( h(\theta )+\pi s^{2}\sum_{i=2}^{n}%
\widehat{\varepsilon }_{i}x_{i}^{2}+R_{3}\left( sx_{2},sx_{3},\ldots
sx_{n}\right) \right) \\
&=&h(\theta )+x_{1}+\pi \sum_{i=2}^{n}\widehat{\varepsilon }%
_{i}x_{i}^{2}+R(x)
\end{eqnarray*}%
where $R(x)=o\left( \left\vert x\right\vert ^{2}\right) $ as $x\rightarrow 0$%
. The would-be factor of $s$ next to $\pi $ can be deleted since the
difference gets absorbed into $R(x)$. What we have just proved is that,%
\begin{equation}
h\left( \theta +x_{1}h(\theta )\theta +Q(x_{2},x_{3},\ldots x_{n})\right)
=h(\theta )+x_{1}+\pi \sum_{i=2}^{n}\widehat{\varepsilon }_{i}x_{i}^{2}+R(x)
\label{local taylor}
\end{equation}%
Since the left side of (\ref{local taylor}) is a $C^{3}$ function of $x$ (in
a neighborhood of $0$), we have in fact that $R(x)=O\left( \left\vert
x\right\vert ^{3}\right) $ as $x\rightarrow 0$. But clearly $R(x)=O\left(
\left\vert x\right\vert ^{2}\right) $ as $x\rightarrow \infty $ and
therefore $\left\vert R(x)\right\vert \leq C\min \left\{ \left\vert
x\right\vert ^{2},\left\vert x\right\vert ^{3}\right\} $ for all $x\in 
\mathbb{R}^{n}$. Note: in the special case where $h$ is a norm and $\mathcal{%
M}$ is everywhere $C^{3}$, it follows that $h$ is $C^{3}$ on $\mathbb{R}%
^{n}\backslash \left\{ 0\right\} $ and the estimate $\left\vert R_{\theta
}\left( x\right) \right\vert \leq C\min \left\{ \left\vert x\right\vert
^{2},\left\vert x\right\vert ^{3}\right\} $ holds for all $\theta \in 
\mathcal{M}$ and $x\in \mathbb{R}^{n}$ with a single $C>0$ independent of $%
\theta $. Using linearity of $Q$ and (\ref{local taylor}),%
\begin{eqnarray*}
&&h\left( t\theta +\frac{x_{1}h(\theta )}{\rho ^{\prime }(th(\theta ))}%
\theta +\sqrt{\frac{t}{\left\vert \rho ^{\prime }(th(\theta ))\right\vert }}%
Q(x_{2},\ldots x_{n})\right) \\
&=&th\left( \theta +\frac{x_{1}h(\theta )}{t\rho ^{\prime }(th(\theta ))}%
\theta +\sqrt{\frac{1}{t\left\vert \rho ^{\prime }(th(\theta ))\right\vert }}%
Q(x_{2},\ldots x_{n})\right) \\
&=&th(\theta )+\frac{x_{1}}{\rho ^{\prime }(th(\theta ))}+\frac{\pi }{%
\left\vert \rho ^{\prime }(th(\theta ))\right\vert }\sum_{i=2}^{n}\widehat{%
\varepsilon }_{i}x_{i}^{2} \\
&&+tR\left( \frac{x_{1}}{t\rho ^{\prime }(th(\theta ))},\sqrt{\frac{1}{%
t\left\vert \rho ^{\prime }(th(\theta ))\right\vert }}x_{2},\ldots ,\sqrt{%
\frac{1}{t\left\vert \rho ^{\prime }(th(\theta ))\right\vert }}x_{n}\right)
\end{eqnarray*}%
Set%
\begin{equation*}
a=h\left( t\theta +\frac{x_{1}h(\theta )}{\rho ^{\prime }(th(\theta ))}%
\theta +\sqrt{\frac{t}{\left\vert \rho ^{\prime }(th(\theta ))\right\vert }}%
Q(x_{2},\ldots x_{n})\right) \hspace{0.65in}b=th(\theta )
\end{equation*}%
By Taylor's theorem there exists $\xi $ between $a$ and $b$ such that%
\begin{eqnarray}
\left\vert \rho (a)-\rho (b)-(a-b)\rho ^{\prime }(b)\right\vert &=&\frac{1}{2%
}\left( a-b\right) ^{2}\left\vert \rho ^{\prime \prime }\left( \xi \right)
\right\vert  \label{random Taylor ee} \\
&=&\frac{1}{2}\left( a-b\right) ^{2}\rho ^{\prime }(b)^{2}\frac{\left\vert
\rho ^{\prime \prime }\left( b+\left( \xi -b\right) \right) \right\vert }{%
\rho ^{\prime }(b)^{2}}
\end{eqnarray}%
Distributing $\rho ^{\prime }(b)$ into the following expression for $a-b$,%
\begin{eqnarray*}
&&\frac{x_{1}}{\rho ^{\prime }(th(\theta ))}+\frac{\pi }{\left\vert \rho
^{\prime }(th(\theta ))\right\vert }\sum_{i=2}^{n}\widehat{\varepsilon }%
_{i}x_{i}^{2} \\
&&+tR\left( \frac{x_{1}}{t\rho ^{\prime }(th(\theta ))},\sqrt{\frac{1}{%
t\left\vert \rho ^{\prime }(th(\theta ))\right\vert }}x_{2},\ldots ,\sqrt{%
\frac{1}{t\left\vert \rho ^{\prime }(th(\theta ))\right\vert }}x_{n}\right)
\end{eqnarray*}%
and applying the triangle inequality, (\ref{random Taylor ee}) implies that%
\begin{equation*}
\left\vert \rho (a)-\rho (b)-x_{1}-\pi \sum_{i=2}^{n}\varepsilon
_{i}x_{i}^{2}\right\vert
\end{equation*}%
is bounded above by%
\begin{eqnarray*}
&&\left\vert \rho ^{\prime }\left( th(\theta )\right) tR\left( \frac{x_{1}}{%
t\rho ^{\prime }(th(\theta ))},\sqrt{\frac{1}{t\left\vert \rho ^{\prime
}(th(\theta ))\right\vert }}x_{2},\ldots ,\sqrt{\frac{1}{t\left\vert \rho
^{\prime }(th(\theta ))\right\vert }}x_{n}\right) \right\vert \\
&&+\frac{1}{2}\left( a-b\right) ^{2}\rho ^{\prime }(b)^{2}\frac{\left\vert
\rho ^{\prime \prime }\left( b+\left( \xi -b\right) \right) \right\vert }{%
\rho ^{\prime }(b)^{2}}
\end{eqnarray*}%
Now $\left\vert \xi -b\right\vert \leq \left\vert a-b\right\vert $, which is
bounded above by%
\begin{eqnarray*}
&&\frac{\pi }{\left\vert \rho ^{\prime }(th(\theta ))\right\vert }\left(
\left\vert x_{1}\right\vert +\sum_{i=2}^{n}x_{i}^{2}\right) \\
&&+\left\vert tR\left( \frac{x_{1}}{t\rho ^{\prime }(th(\theta ))},\sqrt{%
\frac{1}{t\left\vert \rho ^{\prime }(th(\theta ))\right\vert }}x_{2},\ldots ,%
\sqrt{\frac{1}{t\left\vert \rho ^{\prime }(th(\theta ))\right\vert }}%
x_{n}\right) \right\vert
\end{eqnarray*}%
and%
\begin{eqnarray*}
&&\left\vert tR\left( \frac{x_{1}}{t\rho ^{\prime }(th(\theta ))},\sqrt{%
\frac{1}{t\left\vert \rho ^{\prime }(th(\theta ))\right\vert }}x_{2},\ldots ,%
\sqrt{\frac{1}{t\left\vert \rho ^{\prime }(th(\theta ))\right\vert }}%
x_{n}\right) \right\vert \\
&\leq &Ct\min \left\{ \frac{\left\vert x_{1}\right\vert ^{2}}{\left\vert
t\rho ^{\prime }(th(\theta ))\right\vert ^{2}}+\frac{\sum_{i=2}^{n}\left%
\vert x_{i}\right\vert ^{2}}{\left\vert t\rho ^{\prime }(th(\theta
))\right\vert },\frac{\left\vert x_{1}\right\vert ^{3}}{\left\vert t\rho
^{\prime }(th(\theta ))\right\vert ^{3}}+\frac{\left(
\sum_{i=2}^{n}\left\vert x_{i}\right\vert ^{2}\right) ^{3/2}}{\left\vert
t\rho ^{\prime }(th(\theta ))\right\vert ^{3/2}}\right\}
\end{eqnarray*}%
The main estimate then follows from the definition $f(x)=\exp \left( -\rho
\left( h\left( x\right) \right) \right) $. Lastly,%
\begin{eqnarray*}
\left\vert \mathrm{Det}(T_{t})\right\vert &=&\frac{\left\vert \left\langle
\theta ,\nabla h(\theta )\right\rangle \right\vert }{\left\vert \rho
^{\prime }(th(\theta ))\right\vert \cdot \left\vert \nabla h(\theta
)\right\vert }\times \frac{\mathrm{vol}_{n-1}(QB_{2}^{n-1})}{\mathrm{vol}%
_{n-1}(B_{2}^{n-1})}\times \left( \frac{t}{\left\vert \rho ^{\prime
}(th(\theta ))\right\vert }\right) ^{(n-1)/2} \\
&=&\frac{1}{\left\vert \rho ^{\prime }(th(\theta ))\right\vert \cdot
\left\vert \nabla h(\theta )\right\vert }\times \mathrm{Det}(E)\times \left( 
\frac{t}{\left\vert \rho ^{\prime }(th(\theta ))\right\vert }\right)
^{(n-1)/2} \\
&=&\left\vert \kappa (\theta )\right\vert ^{-1/2}\frac{1}{\left\vert \nabla
h(\theta )\right\vert \cdot \left\vert \rho ^{\prime }(th(\theta
))\right\vert }\left( \frac{2\pi t}{\left\vert \nabla h(\theta )\right\vert
\cdot \left\vert \rho ^{\prime }(th(\theta ))\right\vert }\right) ^{(n-1)/2}
\end{eqnarray*}
\end{proof}

\begin{proof}[Proof of Theorem \protect\ref{th mixed homo}]
Note that the function%
\begin{equation*}
h_{1}(x)=\left\{ 
\begin{array}{ccc}
q_{1}(x)\left\vert q_{1}(x)\right\vert ^{-1+1/p(1)} & : & q_{1}(x)\neq 0 \\ 
0 & : & q_{1}(x)=0%
\end{array}%
\right.
\end{equation*}%
has the property that $h_{1}(\alpha x)=\alpha h_{1}(x)$ for all $x\in 
\mathbb{R}^{n}$ and all $\alpha \geq 0$, and%
\begin{equation}
q_{1}(x)=\left\{ 
\begin{array}{ccc}
h_{1}(x)\left\vert h_{1}(x)\right\vert ^{-1+p(1)} & : & h_{1}(x)\neq 0 \\ 
0 & : & h_{1}(x)=0%
\end{array}%
\right.  \label{qqq}
\end{equation}%
Consider the variable $t>0$ large enough so that%
\begin{equation*}
t^{p(1)}>\left\vert \sum_{i=2}^{m}t^{p(i)}q_{i}(\theta )\right\vert
\end{equation*}%
and set%
\begin{eqnarray}
s &=&\left\vert \sum_{i=1}^{m}q_{i}\left( t\theta \right) \right\vert
^{1/p(1)}=(1+o(1))t  \notag \\
q(x) &=&q_{1}(x)+\sum_{i=2}^{m}s^{p(i)-p(1)}q_{i}(x)  \label{defq} \\
\omega &=&ts^{-1}\theta  \notag
\end{eqnarray}%
Then%
\begin{equation*}
q(\omega )=s^{-p(1)}\sum_{i=1}^{m}t^{p(i)}q_{i}(\theta )=\left\vert
\sum_{i=1}^{m}q_{i}\left( t\theta \right) \right\vert
^{-1}\sum_{i=1}^{m}q_{i}\left( t\theta \right) =q_{1}(\theta )
\end{equation*}%
Below we shall refer to functions $R_{i}(\cdot ,\cdot )$. The first variable
may be $z\in \mathbb{R}^{n-1}$ or $x\in \mathbb{R}^{n}$ (this will always be
clear from the context). Each has the property that $\forall t>t_{0}$,%
\begin{equation*}
\lim_{z\rightarrow 0}\left\vert z\right\vert ^{-2}R_{i}(z,t)=0
\end{equation*}%
i.e. $R_{i}(z,t)=o\left( \left\vert z\right\vert ^{2}\right) $ as $%
z\rightarrow 0$. However the rate of convergence may (possibly) depend on $t$%
. We shall also refer to functions $\delta $, $\delta _{i}$, $\delta _{i,j}$%
, $\gamma _{i}$ and $\gamma _{i,j}$. These quantities (coefficients) are
functions of $t$ that do not depend on $z\in \mathbb{R}^{n-1}$ or $x\in 
\mathbb{R}^{n}$, such that as $t\rightarrow \infty $ we have $\delta
\rightarrow 0$ and $\gamma \rightarrow 1$. They may also denote different
functions from one appearance to the next. Let $(v_{i})_{1}^{n-1}$ denote
the principal directions of $\mathcal{M}$ at $\theta $, and $(\kappa
_{i})_{1}^{n-1}$ the corresponding principal curvatures. Let $W:\mathbb{R}%
^{n-1}\rightarrow \nabla q_{1}(\theta )^{\perp }$ be the linear isometry
such that $We_{i}=v_{i}$ ($1\leq i\leq n-1$). Let $U\in SO(n)$ such that%
\begin{equation*}
U\left( \frac{\nabla q_{1}(\theta )}{\left\vert \nabla q_{1}(\theta
)\right\vert }\right) =\frac{\nabla q(\omega )}{\left\vert \nabla q(\omega
)\right\vert }
\end{equation*}%
and such that $Ux=x$ for all $x\in \left\{ \nabla q_{1}(\theta ),\nabla
q(\omega )\right\} ^{\bot }$. As $t\rightarrow \infty $, $\omega \rightarrow
\theta $ and since $q_{1}$ is $C^{1}$ at $\theta $, $\nabla q_{1}(\omega
)\rightarrow \nabla q_{1}(\theta )$. By (\ref{defq}) and the fact that $%
q_{i} $ ($1\leq i\leq m$) are $C^{1}$ at $\theta $ and do not depend on $t$,
it follows that $\nabla q(\omega )\rightarrow \nabla q_{1}(\theta )$, and
therefore if $t$ is sufficiently large $\nabla q(\omega )\neq 0$ and $U$ is
indeed well defined. Furthermore, $U\rightarrow I_{n}$ (the identity matrix)
in the standard topology on $\mathbb{R}^{n\times n}$. It follows as in the
proof of Theorem \ref{th homo level}, see in particular (\ref{some approxe}%
), that%
\begin{equation}
h_{1}\left( \theta +W(z)\right) =h_{1}(\theta )+\frac{1}{2}h_{1}(\theta
)\left\vert \nabla h_{1}(\theta )\right\vert \sum_{i=1}^{n-1}\kappa
_{i}z_{i}^{2}+R_{1}(z)  \label{nomorenames}
\end{equation}%
Since $(1+\varepsilon )^{p}=1+p\varepsilon +o(\varepsilon )$ as $\varepsilon
\rightarrow 0$, (\ref{nomorenames}) implies that%
\begin{equation*}
q_{1}\left( \theta +W(z)\right) =q_{1}(\theta )+\sum_{i=1}^{n-1}\beta
_{i}z_{i}^{2}+R_{2}(z,t)
\end{equation*}%
where $\beta _{i}\neq 0$ ($1\leq i\leq n-1$). Since each $q_{i}$ is $C^{2}$
at $\theta $, $H_{q}(\omega )\rightarrow H_{q_{1}}(\theta )$. Using this and
the fact that $\left\vert x-U(x)\right\vert <\delta (t)\left\vert
x\right\vert $ for an appropriate $\delta (\cdot )$, and the fact that $%
Range(UW)=\nabla q(\omega )^{\bot }$,%
\begin{equation*}
q\left( \omega +UW(z)\right) =q(\omega )+\sum_{i=1}^{n-1}\beta _{i}\gamma
_{i}z_{i}^{2}+\sum_{1\leq i<j\leq n-1}\delta _{i,j}z_{i}z_{j}+\widehat{R}%
_{2}(z,t)
\end{equation*}%
In particular, the Hessian matrix of the function $z\mapsto q\left( \omega
+UW(z)\right) $ is invertible when evaluated at $z=0$, and there exists an
injective linear map $Q:\mathbb{R}^{n-1}\rightarrow \nabla q(\omega )^{\perp
}$, such that%
\begin{equation}
q\left( \omega +Q(z)\right) =q(\omega )+\pi \sum_{i=1}^{n-1}\varepsilon
_{i+1}z_{i}^{2}+R_{3}(z,t)  \label{co3}
\end{equation}%
The vector $(\varepsilon _{i})_{2}^{n}\in \left\{ \pm 1\right\} ^{n-1}$ is
indexed by $\{2,\ldots ,n\}$, hence the need for the subscript $i+1$ in the
summation. The second order Taylor expansion in $n$ variables can be written
as%
\begin{eqnarray*}
&&q\left( \omega +x_{1}\omega +Q(x_{2},\ldots ,x_{n})\right) \\
&=&q(\omega )+\left\langle \nabla q(\omega ),\omega \right\rangle x_{1}+%
\frac{1}{2}x_{1}^{2}\left\langle \omega ,H_{q}(\omega )\omega \right\rangle
+x_{1}\left\langle \omega ,H_{q}(\omega )Q(x_{2},\ldots ,x_{n})\right\rangle
\\
&&+\frac{1}{2}\left\langle Q(x_{2},\ldots ,x_{n}),H_{q}(\omega
)Q(x_{2},\ldots ,x_{n})\right\rangle +R_{5}(x,t)
\end{eqnarray*}%
Write $\left\langle \omega ,H_{q}(\omega )Q(x_{2},\ldots
,x_{n})\right\rangle =\sum_{2}^{n}a_{i}x_{i}$. $Q$ can be chosen so that as $%
t\rightarrow \infty $, $Q=Q(t)$ converges in $\mathbb{R}^{n\times (n-1)}$
and therefore so does each $a_{i}$. Combining this with (\ref{co3}) gives%
\begin{eqnarray*}
&&q\left( \omega +x_{1}\omega +Q(x_{2},\ldots ,x_{n})\right) \\
&=&q(\omega )+\left\langle \nabla q(\omega ),\omega \right\rangle x_{1}+\pi
\sum_{i=2}^{n}\varepsilon _{i}x_{i}^{2}+\frac{1}{2}\left\langle \omega
,H_{q}(\omega )\omega \right\rangle
x_{1}^{2}+x_{1}\sum_{i=2}^{n}a_{i}x_{i}+R_{4}(x,t)
\end{eqnarray*}%
Using Taylor's theorem, the $C^{3}$ condition, and (\ref{defq}), it follows
that we may take $\left\vert R_{4}(x,t)\right\vert \leq c\left\vert
x\right\vert ^{3}$ provided $\left\vert x\right\vert <\varepsilon $, where $%
c,\varepsilon >0$ do not depend on $t$. Note that $\left\langle \omega
,\nabla q(\omega )\right\rangle \rightarrow \left\langle \theta ,\nabla
q_{1}(\theta )\right\rangle \neq 0$. Using linearity of $Q$ and recalling
that $t\theta =s\omega $,%
\begin{eqnarray*}
&&\sum_{i=1}^{m}q_{i}\left( t\theta +\frac{x_{1}}{s^{p(1)-1}\left\langle
\omega ,\nabla q(\omega )\right\rangle }\omega +\frac{1}{s^{p(1)/2-1}}%
Q(x_{2},\ldots x_{n})\right) \\
&=&s^{p(1)}q\left( \omega +\frac{x_{1}}{s^{p(1)}\left\langle \omega ,\nabla
q(\omega )\right\rangle }\omega +\frac{1}{s^{p(1)/2}}Q(x_{2},\ldots
x_{n})\right) \\
&=&s^{p(1)}\left[ q(\omega )+\frac{x_{1}}{s^{p(1)}}+\pi \sum_{i=2}^{n}\frac{%
\varepsilon _{i}x_{i}^{2}}{s^{p(1)}}+\frac{\left\langle \omega ,H_{q}(\omega
)\omega \right\rangle x_{1}^{2}}{2s^{2p(1)}\left\langle \omega ,\nabla
q(\omega )\right\rangle ^{2}}+\frac{x_{1}\sum_{i=2}^{n}a_{i}x_{i}}{%
s^{3p(1)/2}\left\langle \omega ,\nabla q(\omega )\right\rangle }\right] \\
&&+s^{p(1)}R_{4}\left( \frac{x_{1}}{s^{p(1)}\left\langle \omega ,\nabla
q(\omega )\right\rangle },s^{-p(1)/2}x_{2},s^{-p(1)/2}x_{3},\ldots
,s^{-p(1)/2}x_{n},t\right) \\
&=&\sum_{i=1}^{m}q_{i}\left( t\theta \right) +x_{1}+\pi
\sum_{i=2}^{n}\varepsilon _{i}x_{i}^{2}+\frac{\left\langle \omega
,H_{q}(\omega )\omega \right\rangle x_{1}^{2}}{2s^{p(1)}\left\langle \omega
,\nabla q(\omega )\right\rangle ^{2}}+\frac{x_{1}\sum_{i=2}^{n}a_{i}x_{i}}{%
s^{p(1)/2}\left\langle \omega ,\nabla q(\omega )\right\rangle } \\
&&+s^{p(1)}R_{4}\left( \frac{x_{1}}{s^{p(1)}\left\langle \omega ,\nabla
q(\omega )\right\rangle },s^{-p(1)/2}x_{2},s^{-p(1)/2}x_{3},\ldots
,s^{-p(1)/2}x_{n},t\right)
\end{eqnarray*}%
The result then follows from the definition of $f$.
\end{proof}

\begin{proof}[Proof of Theorem \protect\ref{find scale}]
Consider the bilinear form on $\mathbb{R}^{n}$ defined by%
\begin{equation*}
\left\langle x,z\right\rangle _{f,y}=\frac{1}{2\pi }\left\langle
x,H_{g}(y)z\right\rangle +\left\langle \nabla g(y),x\right\rangle \cdot
\left\langle \nabla g(y),z\right\rangle =\left\langle x,Az\right\rangle
=\left\langle Bx,Bz\right\rangle
\end{equation*}%
where $g=-\log f$ and $B=A^{1/2}$ (in a moment we will see that $A$ is
positive definite). This bilinear form has two key properties. The first
property is that it does not depend on the underlying coordinate structure
of $\mathbb{R}^{n}$: If $W:\mathbb{R}^{n}\rightarrow \mathbb{R}^{n}$ is any
linear bijection and $u\in \mathbb{R}^{n}$ then%
\begin{equation*}
\left\langle x,z\right\rangle _{f,y}=\left\langle Wx,Wz\right\rangle _{%
\widetilde{f},u+Wy}
\end{equation*}%
where $\widetilde{f}(x)=f\left( W^{-1}(x-u)\right) $. This follows since%
\begin{eqnarray*}
\nabla \widetilde{g}(x) &=&\left( W^{-1}\right) ^{T}\nabla g\left(
W^{-1}(x-u)\right) \\
H_{\widetilde{g}}(x) &=&\left( W^{-1}\right) ^{T}H_{g}\left(
W^{-1}(x-u)\right) W^{-1}
\end{eqnarray*}%
where $\widetilde{g}(x)=g\left( W^{-1}(x-u)\right) $. The second is that%
\begin{equation*}
\left\langle \cdot ,\cdot \right\rangle _{\exp (-\xi ),0}=\left\langle \cdot
,\cdot \right\rangle
\end{equation*}%
i.e. when $f(x)=\exp (-\xi (x))$ and $y=0$, it reduces to the standard
Euclidean inner product. Therefore, setting $W=T^{-1}$ and $u=-Wy$,%
\begin{equation*}
\left\langle x,z\right\rangle _{f,y}=\left\langle
T^{-1}x,T^{-1}z\right\rangle _{\exp (-\xi ),0}=\left\langle
T^{-1}x,T^{-1}z\right\rangle
\end{equation*}%
Thus $\left\langle \cdot ,\cdot \right\rangle _{f,y}$ is an inner product
and $A$ is a symmetric positive definite matrix. This implies that $A$ and $%
B $ are invertible. Since $\left\langle Bx,Bz\right\rangle =$ $\left\langle
T^{-1}x,T^{-1}z\right\rangle $ for all $x,z\in \mathbb{R}^{n}$, $T=B^{-1}U$
for some $U\in O(n)$. Since (\ref{ph}) can be written as%
\begin{equation*}
g(y+Tx)-g(y)=\xi (x)
\end{equation*}%
it follows that%
\begin{equation}
\nabla _{x}g(y+Tx)=\nabla \xi (x)  \label{qa}
\end{equation}%
where $\nabla _{x}g(y+Tx)$ denotes the gradient of the function $x\mapsto
g(y+Tx)$. However by the chain rule%
\begin{equation}
\nabla _{x}g(y+Tx)=T^{T}\nabla g(y+Tx)  \label{qb}
\end{equation}%
Setting $x=0$ and equating the right sides of (\ref{qa}) and (\ref{qb}), $%
e_{1}=U^{T}B^{-1}\nabla g(y)$ and the first result follows by setting $%
F=U^{T}$. We now consider the special case $\xi (x)=x_{1}+\pi
\sum_{2}^{n}x_{i}^{2}$ and a matrix $G\in O(n)$ such that $GA^{-1/2}\nabla
g(y)=e_{1}$. Now $FG^{T}\in O(n)$ and by (\ref{ro1}) and (\ref{ro2}),%
\begin{equation*}
FG^{T}e_{1}=FG^{T}GA^{-1/2}\nabla g(y)=FA^{-1/2}\nabla g(y)=e_{1}
\end{equation*}%
It then follows by rotational invariance of $\Lambda _{n}$ about the $x_{1}$
direction that%
\begin{equation*}
\frac{f\left( y+A^{-1/2}G^{T}x\right) }{f(y)}=\frac{f\left(
y+A^{-1/2}F^{T}FG^{T}x\right) }{f(y)}=\Lambda _{n}\left( FG^{T}x\right)
=\Lambda _{n}\left( x\right)
\end{equation*}
\end{proof}

\begin{proof}[Proof of Theorem \protect\ref{self sim}]
Consider the Hilbert space $\mathfrak{H}=\nabla \Lambda (y)^{\bot }=\left(
\nabla \log \Lambda (y)\right) ^{\bot }$ endowed with the inner product and
corresponding norm%
\begin{equation*}
\left\langle x,y\right\rangle _{\mathfrak{H}}=\sum_{i=2}^{n}x_{i}y_{i}%
\hspace{0.65in}\left\Vert x\right\Vert _{\mathfrak{H}}=\left(
\sum_{i=2}^{n}x_{i}^{2}\right) ^{1/2}
\end{equation*}%
The function $\left\Vert \cdot \right\Vert _{\mathfrak{H}}$ is indeed a norm
because if $\left\Vert x\right\Vert _{\mathfrak{H}}=0$ then $x=te_{1}$ for
some $t\in \mathbb{R}$, and it follows from the definition of $\mathfrak{H}$
that $t=0$. Consider a linear map $U:\mathbb{R}^{n-1}\rightarrow \mathfrak{H}
$ such that for all $x\in \mathbb{R}^{n-1}$, $\left\Vert Ux\right\Vert _{%
\mathfrak{H}}=|x|$. Defining $Tx=x_{1}e_{1}+U(x_{2},\ldots x_{n})$, we then
have%
\begin{eqnarray*}
\log \frac{\Lambda _{n}\left( y+Tx\right) }{\Lambda _{n}(y)} &=&-x_{1}-\pi
\sum_{i=2}^{n}x_{i}^{2}-e_{1}^{\sharp }\left( U(x_{2},\ldots ,x_{n})\right)
\\
&&-2\pi \sum_{i=2}^{n}y_{i}e_{i}^{\sharp }\left( U(x_{2},\ldots
,x_{n})\right) \\
&=&\log \Lambda _{n}(x)+\left\langle U(x_{2},\ldots ,x_{n}),\nabla \log
\Lambda _{n}(y)\right\rangle =\log \Lambda _{n}(x)
\end{eqnarray*}%
where $e_{i}^{\sharp }(z)=z_{i}$ is the $i^{th}$ coordinate functional.
\end{proof}

\section*{Acknowledgements}

Many thanks to Alexander Koldobsky for providing a reference on intersection
bodies and to Elisabeth Werner for asking a question about affine surface
area for measures that was part of our initial motivation for studying the
tail behavior of functions. Parts of this paper were written while the
author was a PhD student at the University of Missouri, a Gibbs Assistant
Professor at Yale University, and a post-doctoral fellow at the Weizmann
Institute of Science.

\end{document}